\documentclass[11pt,a4paper]{article}
\usepackage[T1]{fontenc}
\usepackage{lmodern}
\usepackage[a4paper,left=0.94in,right=0.94in,top=0.80in,bottom=0.79in]{geometry}
\usepackage{amsmath,amssymb,amsthm,mathtools}
\usepackage{microtype}
\usepackage{booktabs,array}
\usepackage{fancyvrb}
\usepackage{upquote}
\usepackage{enumitem}
\usepackage[numbers,sort&compress]{natbib}
\usepackage{needspace}
\usepackage{float}
\usepackage{etoolbox}
\usepackage{hyperref}
\hypersetup{hidelinks,pdftitle={The cofinite regularity spectrum of F-irregular graphs},pdfauthor={Zhanhe Zhang}}
\usepackage{url}
\BeforeBeginEnvironment{theorem}{\Needspace{7\baselineskip}}
\BeforeBeginEnvironment{lemma}{\Needspace{6\baselineskip}}
\BeforeBeginEnvironment{corollary}{\Needspace{5\baselineskip}}
\numberwithin{equation}{section}
\newtheorem{theorem}{Theorem}[section]
\newtheorem{lemma}[theorem]{Lemma}
\newtheorem{corollary}[theorem]{Corollary}

\newcommand{\FF}{\mathbb F}
\newcommand{\RR}{\mathbb R}
\newcommand{\Hom}{\operatorname{hom}}
\newcommand{\Inj}{\operatorname{inj}}
\newcommand{\Aut}{\operatorname{Aut}}
\newcommand{\Spec}{\mathcal R}
\newcommand{\diag}{\operatorname{diag}}
\newcommand{\eps}{\varepsilon}
\newcommand{\one}{\boldsymbol1}
\newcommand{\istar}{i_*}

\selectfont
\newif\ifprintseedcode
\printseedcodefalse 
\begin{document}
\begin{center}
  {\fontsize{17.2}{20}\selectfont\bfseries The cofinite regularity spectrum of $F$-irregular graphs\par}
  \vspace{1.1em}
  {\normalsize Zhanhe Zhang\footnotemark\par}
  \vspace{0.55em}
  {\small School of International Trade and Economics,\par}
  {\small Central University of Finance and Economics, Beijing, China\par}
  \vspace{0.95em}
  {\normalsize October 2026\par}
\end{center}
\footnotetext{Corresponding author. Email: \texttt{2024311337@email.cufe.edu.cn}.}
\vspace{0.1em}
\begin{center}\textbf{Abstract}\end{center}
\vspace{-0.55em}
\noindent
For a connected graph $F$, a graph $G$ is $F$-irregular if the numbers of copies of $F$ containing its vertices are pairwise distinct. We characterize which connected patterns admit regular $F$-irregular hosts of every sufficiently large degree. Stars admit no regular $F$-irregular host; for every other $F$, every sufficiently large integer $r$ is realized by a connected $r$-regular $F$-irregular graph of diameter two and order $r+o(r)$, optimal to first order. For triangles, and equivalently four-vertex paths in regular graphs, we construct such hosts of order less than $r+73r^{\log29/\log58}$ above an explicit threshold, with the constant reduced to $44$ for even $r$. Lexicographic words encode distinct rooted subgraph counts, while layered connections allow the common degree to be prescribed. Rooted-homomorphism expansions and spectral amplification establish separation for general patterns; an exact uniform-increment identity gives the quantitative triangle construction. The finite seed graphs are explicitly constructed, and their required rooted-count distinctions are verified by exact-arithmetic certificates.

\vspace{0.9em}
\par\vspace{0.35em}
\par\vspace{1.1em}

\section{Introduction}\label{sec:intro}
Vertices of a regular graph have equal degree, yet their rooted counts of a fixed graph $F$ need not agree. For graphs $F,G$ and $v\in V(G)$, write
\[
 d_F^G(v)=\#\{H\subseteq G:H\cong F,\ v\in V(H)\},
\]
where copies need not be induced. We call $G$ \emph{$F$-irregular} if the values $d_F^G(v)$ are pairwise distinct. The questions are then which common degrees can be prescribed and, for a prescribed degree $r$, how close the order of the host can come to the elementary lower bound $r+1$.

Stars give an immediate obstruction. In an $r$-regular graph,
\begin{equation}\label{eq:star}
 d_{K_{1,s}}^G(v)=(s+1)\binom rs\qquad(s\ge2)
\end{equation}
is independent of $v$. To exclude vacuous hosts, require $|V(G)|\ge|V(F)|$ and define the \emph{regularity spectrum}
\begin{equation}\label{eq:spectrum}
 \Spec(F)=\left\{r\in\mathbb Z_{\ge0}:
 \begin{array}{l}
 \text{there is a connected $r$-regular $F$-irregular graph $G$}\\
 \text{with $|V(G)|\ge|V(F)|$}
 \end{array}\right\}.
\end{equation}
For $r\in\Spec(F)$ let $n_F(r)$ be the minimum order of a graph witnessing this membership, and put $n_F(r)=+\infty$ otherwise. The elementary lower bound $n_F(r)\ge r+1$ makes order minimization a natural companion to the regularity-spectrum problem.

\begin{theorem}[Cofinite regularity spectrum]\label{thm:main}
 Let $F$ be a finite, simple, connected graph with at least three vertices.  If $F$ is a star, then $\Spec(F)=\varnothing$.  Otherwise there exists an integer $r_F$ such that for \emph{every} integer $r\ge r_F$ there is a finite simple connected $r$-regular $F$-irregular graph $G_r$ of diameter exactly two.  Moreover, the hosts can be chosen so that $|V(G_r)|/r\to1$ as $r\to\infty$.
\end{theorem}

\begin{corollary}[Existence classification]\label{cor:classification}
 A finite connected graph $F$ of order at least three admits a connected regular $F$-irregular host if and only if $F$ is not a star.  In the non-star case there are infinitely many such hosts of odd regularity and infinitely many of even regularity, each of diameter two.
\end{corollary}
\begin{proof}
 Theorem~\ref{thm:main} gives all sufficiently large integers, including both parities, and the star obstruction is~\eqref{eq:star}.  Graphs of different regularities are nonisomorphic.
\end{proof}

Theorem~\ref{thm:main} also settles the first-order extremal problem. For every fixed connected non-star $F$ of order at least three,
\begin{equation}\label{eq:extremal-first-order}
 n_F(r)=r+o(r),\qquad \lim_{r\to\infty}\frac{n_F(r)}r=1.
\end{equation}
This is optimal to first order because every simple $r$-regular graph has at least $r+1$ vertices. The conclusion holds while the host is required to be connected and have diameter two.

For triangles, the surplus above this elementary bound can be controlled by an explicit power of $r$. Write $\tau_G(v)=d_{K_3}^G(v)$; a graph with injective $\tau_G$ is also called \emph{triangle-distinct}. In an $r$-regular graph,
\[
 d_{P_4}^G(v)=2r(r-1)^2-4\tau_G(v),
\]
so $K_3$-irregularity and $P_4$-irregularity are equivalent. The identity is proved in Lemma~\ref{lem:basic-identities}. Set
\begin{equation}\label{eq:qt-data}
 \beta_\triangle=\frac{\log29}{\log58}=0.829292720\ldots,
 \qquad r_\triangle=58^{11}-\frac{29^{11}+1}{2}+1.
\end{equation}
\begin{theorem}[Quantitative order bound]\label{thm:quantitative}
For every integer $r\ge r_\triangle$, there is a finite simple connected $r$-regular graph $G_r$ of diameter exactly two such that both $d_{K_3}^{G_r}$ and $d_{P_4}^{G_r}$ are injective and
\begin{equation}\label{eq:qt-bound}
 r+1\le |V(G_r)|<r+73r^{\beta_\triangle}.
\end{equation}
For even $r$, the same requirements can be met with $|V(G_r)|<r+44r^{\beta_\triangle}$. The construction has an explicit adjacency rule.
\end{theorem}
For triangles, the known degree obstruction, recalled in~\eqref{eq:qt-lower}, and Theorem~\ref{thm:quantitative} give, for all sufficiently large $r$,
\[
 3+\sqrt{2r+4}\le n_{K_3}(r)-r<73r^{\beta_\triangle}.
\]
This comparison places the quantitative result in an extremal-order problem: its constructive exponent is strictly below one but remains above the known square-root obstruction. The power bound concerns $K_3$ and, by the regular-graph identity, $P_4$; for other connected non-star patterns we obtain the first-order estimate~\eqref{eq:extremal-first-order}. The more general seed-lifting principle, Theorem~\ref{thm:qt-lifting}, yields exponents $\log q/\log(aq)$ from $a$ equal-degree seeds of order $q$, while Theorem~\ref{thm:qt-alphabet} restricts this class of inputs by $q\ge8a+8$.

\subsection{Related work}
The $F$-degree notion was introduced by Chartrand, Holbert, Oellermann and Swart~\cite{CHOS1987}; irregularity questions in graph theory were subsequently formulated by Chartrand, Erd\H{o}s and Oellermann~\cite{CEO1988}. Draganova~\cite{Draganova2009} studied $F$-continuity, a different local requirement that adjacent vertices have $F$-degrees differing by at most one. Her regular-host counterexamples to $F$-continuity do not address pairwise separation across all vertices.

Without fixing the host degree, Schreib~\cite{Schreib2026} constructed $F$-irregular graphs for connected noncomplete patterns, with complete patterns covered by earlier results. Dovzhenok and Filuta~\cite{DF2026} obtained unrestricted $F$-irregular hosts for all connected patterns of order at least three; Dovzhenok~\cite{Dovzhenok2026} studied diameter-three constructions for patterns of diameter two. Prescribing the ordinary degree adds an independent global constraint: the same $r$ must be realized at every vertex while the rooted $F$-counts remain distinct.

For triangles, Berikkyzy et al.~\cite{Berikkyzy2024} raised the regular-existence question, while Stevanovi\'c et al.~\cite{Stevanovic2024} and Hak, Kozerenko and Serdiuk~\cite{Hak2025} gave regular examples. Zhang~\cite{ZhangK3} constructed regular triangle-distinct graphs for every $r\ge9$, using $2r+O(1)$ vertices for all sufficiently large $r$; together with the nonexistence results of Hak, Kozerenko and Serdiuk~\cite{Hak2025,HakNo8}, this establishes the exact threshold $9$. The threshold-switch method of~\cite{ZhangK3} supplies the $125$-vertex seeds used in our general-pattern construction. A separate complementary $29$-vertex pair yields the sublinear order surplus in Theorem~\ref{thm:quantitative}.

\subsection{Proof strategy}
We construct the host as the complement of an $\omega$-regular graph $R$ on $N$ vertices, with
\[
 r=N-1-\omega.
\]
When $\omega=o(N)$, this immediately yields $N=r+o(r)$; the complement also has diameter two for all sufficiently large $r$. The difficulty is to prescribe $\omega$ exactly while separating every pair of rooted $F$-counts in one auxiliary graph.

The construction has three structural ingredients. First, lexicographic words encode distinct rooted statistics through their highest differing letters. Second, nested suffix-shell connections interpolate the common degree while adding the same constant to every rooted triangle count (Lemma~\ref{lem:shell-invariance}). Third, the rooted complement expansion identifies the cyclic statistics governing a general pattern. Block centering removes rooted tree contributions, and spectral amplification controls the remaining shortest-cycle term. Only in the large-girth case is a further scale-preserving change needed to avoid cancellation of two leading terms. Table~\ref{tab:leading-patterns} records the resulting cases.

Sections~\ref{sec:prelims}--\ref{sec:proof} prove the classification theorem; Sections~\ref{sec:quantitative}--\ref{sec:alphabet} establish the quantitative triangle theorem and its seed-lifting generalization. The appendices give the finite constructions and the additional counting identities. All choices are made before comparing any pair of roots, as formalized in Section~\ref{sec:proof}.

\begin{table}[H]
\centering
\caption{The leading rooted statistics in the four cases of the classification theorem.}
\label{tab:leading-patterns}
\small
\begin{tabular}{@{}p{0.30\linewidth}p{0.24\linewidth}p{0.35\linewidth}@{}}
\toprule
Pattern $F$ & Statistic in $R$ & Principal ingredient \\
\midrule
Contains a triangle & $\tau_R$ & An uncontracted triangle support \\
Non-star tree & $\tau_R$ & Identifying the endpoints of a $P_4$ \\
Triangle-free, contains $C_4$ & $c_4^R$ & A bipartite seed alphabet \\
Girth $g\ge5$ & $\tau_R$ and $z_g^R$ & Amplification; phase selection for odd $g$ \\
\bottomrule
\end{tabular}
\end{table}

\section{Rooted counts and finite seed alphabets}\label{sec:prelims}
All graphs are finite, simple and undirected, and logarithms are natural unless their base is indicated. Write $[n]=\{1,\ldots,n\}$. When comparing roots $v,w$ in the same graph, we use $\Delta f=f(v)-f(w)$; when the graphs differ, the two rooted quantities will be specified. For a rooted graph $(H,a)$, let $\Hom_v(H,a;R)$ and $\Inj_v(H,a;R)$ denote the numbers of homomorphisms and injective homomorphisms, respectively, sending $a$ to $v$.  Put $a_F=|\Aut(F)|$.  Thus
\begin{equation}\label{eq:rootedinj}
 a_Fd_F^G(v)=\sum_{a\in V(F)}\Inj_v(F,a;G).
\end{equation}
For an $\omega$-regular graph $R$ write $\tau_R(v)=d_{K_3}^R(v)$, $c_4^R(v)=d_{C_4}^R(v)$, and
\begin{equation}\label{eq:walk}
 z_s^R(v)=(A_R^s)_{vv}.
\end{equation}
Then $z_3^R=2\tau_R$.  Write $p_4(F)$ for the number of (not necessarily induced) copies of the four-vertex path $P_4$ in $F$, and $c_s(F)$ for the number of cycles $C_s$.

If $X$ and $Y$ are regular graphs of degrees $d_X$ and $d_Y$, respectively, and $|V(Y)|=q$, then their lexicographic product $X[Y]$ satisfies
\begin{equation}\label{eq:lex}
 |V(X[Y])|=|V(X)|q,\qquad d(X[Y])=qd_X+d_Y.
\end{equation}
\subsection{Three counting identities}
The following identities are used by both constructions. The triangle identities are standard; see also the formulas and attributions in~\cite{Berikkyzy2024}. We include the regular-host calculations to fix the counting convention.
\begin{lemma}[Regular-host identities]\label{lem:basic-identities}
Let $H$ be a $d$-regular graph of order $N$. Then
\begin{align}
 \tau_{\overline H}(v)
 &=\binom{N-1}{2}-\frac32d(N-1-d)-\tau_H(v),\label{eq:triangle-complement}\\
 d_{P_4}^H(v)&=2d(d-1)^2-4\tau_H(v).\label{eq:p4exact}
\end{align}
If $X,Y$ are regular and $|V(Y)|=q$, then
\begin{equation}\label{eq:trianglelex}
 \tau_{X[Y]}(x,y)=q^2\tau_X(x)+\tau_Y(y)+\tfrac32d_Xd_Yq.
\end{equation}
\end{lemma}
\begin{proof}
Put $S=N_H(v)$ and $T=V(H)\setminus(S\cup\{v\})$. Summing degrees on $S$ gives $e_H(S,T)=d(d-1)-2\tau_H(v)$, and hence
\[
 e_H(T)=\frac{Nd}{2}-d^2+\tau_H(v).
\]
Since the neighbors of $v$ in $\overline H$ are exactly $T$, subtracting this expression from $\binom{N-1-d}{2}$ gives~\eqref{eq:triangle-complement}.

There are $d(d-1)^2$ oriented non-backtracking walks of length three with $v$ in any prescribed position. A noninjective one closes a triangle through $v$, and there are $2\tau_H(v)$ such walks. Summing over the four positions and dividing by the two orientations of a path proves~\eqref{eq:p4exact}.

In $X[Y]$, a triangle through $(x,y)$ lies in one fiber, meets exactly two fibers, or meets three fibers. These contributions are respectively $\tau_Y(y)$, $d_Xe(Y)+d_Xd_Yq$, and $q^2\tau_X(x)$. Since $e(Y)=qd_Y/2$, their sum proves~\eqref{eq:trianglelex}.
\end{proof}
In particular, complementation preserves triangle-distinctness of regular graphs, and $K_3$- and $P_4$-irregularity are equivalent on this class.

\subsection{The seeds for the general-pattern construction}
Three features of the finite seeds enter the proof. Equal-degree triangle seeds with jointly distinct rooted counts supply the letters used for root separation. For triangle-free patterns containing a four-cycle, a bipartite pair supplies jointly distinct rooted four-cycle counts while its triangle counts vanish. The padding alphabet also has to admit a change of internal scales at fixed total order: here its size $b=31250$ equals $2n^2$ for the upper alphabet size $n=125$. The bipartite pair is obtained from the same complementary triangle seeds by a deterministic $5\times5$ block substitution. Lemma~\ref{lem:seeds} records the finite data; Appendix~\ref{app:embedded} gives their reconstruction rules and exact checks.

\begin{lemma}[Seed certificates]\label{lem:seeds}
 There exist graphs satisfying the following properties.
 \begin{enumerate}[label=(\roman*)]
 \item Two complementary $125$-vertex $62$-regular graphs $A,B$ whose combined $250$ rooted triangle counts are pairwise distinct, with total span $932$.
 \item The $31250$-vertex $7812$-regular graph
 \begin{equation}\label{eq:padseed} H_0=A[A]\ \dot\cup\ B[B] \end{equation}
 has $31250$ pairwise distinct rooted triangle counts, with total span at most $14\,563\,432$.
 \item Two $1250$-vertex $310$-regular bipartite graphs $U,V$ with disjoint edge sets under a common vertex labeling. Their $2500$ rooted four-cycle counts are pairwise distinct and have total span $834\,684$.
 \end{enumerate}
\end{lemma}
\begin{proof}
Appendix~\ref{app:embedded} gives a short explicit construction of the finite graphs and a self-contained exact-arithmetic check of (i) and (iii). Regularity and edge-disjointness follow structurally from the construction; the rooted-count injectivity is also checked there. Part (ii) follows from the lexicographic-product identity~\eqref{eq:trianglelex}.
Consequently, for $S\in\{A,B\}$ the triangle count at $(x,y)\in S[S]$ equals
\[
 15625\tau_S(x)+\tau_S(y)+720750.
\]
Suppose two such values (possibly taken from different components $S[S]$ and $T[T]$) were equal.  Then
\[
 15625\bigl(\tau_S(x)-\tau_T(x')\bigr)
 =\tau_T(y')-\tau_S(y).
\]
The absolute value on the right is at most $932<15625$.  Therefore $\tau_S(x)=\tau_T(x')$, and by the $250$-value injectivity in (i) this forces $S=T$ and $x=x'$.  Equality of the low digits then forces $y=y'$.  Finally the total span is at most $(15625+1)\cdot932=14\,563\,432$, and the degree is $125\cdot62+62=7812$.  This proves (ii).
\end{proof}

We use the pair $(A,B)$ for upper triangle digits, $H_0$ for lower digits, and the pair $(U,V)$ for four-cycle separation. The equal seed degrees keep the word blocks regular.

\section{Layered word graphs and degree interpolation}\label{sec:wiring}
We now place the seed alphabets of Section~\ref{sec:prelims} in a common lexicographic word construction and define the suffix kernels used to connect its blocks. The radices may vary within a fixed finite range.

Let $n_0,\ldots,n_{P-1}$ be integer radices satisfying $2\le n_i\le n_{\max}$ for a fixed $n_{\max}$, ordered from least to most significant, and set
\[
 q_0=1,\quad q_i=\prod_{a<i}n_a,\quad Q=\prod_{a<P}n_a=q_P.
\]
Let $\Gamma=\FF_2^\ell$, whose coordinates are label bits indexed by $0,\ldots,\ell-1$; each $e\in\Gamma$ labels one word block. At each word position $i$, fix two $d_i$-regular graphs $S_{i,0},S_{i,1}$ on the common labeled vertex set $[n_i]$, and a linear map $\chi_i:\Gamma\to\FF_2$; the zero map is allowed. For $e\in\Gamma$, let
\[
 W_e=S_{P-1,\chi_{P-1}(e)}[\cdots[S_{0,\chi_0(e)}]\cdots].
\]
Define the suffix degrees and the common block degree by
\begin{equation}\label{eq:worddegree}
 D_i=\sum_{a<i}d_aq_a,\qquad D=D_P=\sum_{i=0}^{P-1}d_iq_i.
\end{equation}
A vertex has coordinates $(e,x_{P-1},\ldots,x_0)$. To distinguish the word-block labels, we require the joint color map to be injective:
\[
 e\longmapsto\bigl(\chi_0(e),\ldots,\chi_{P-1}(e)\bigr).
\]
The choices below satisfy this condition. For different roots $v=(e,x)$ and $w=(f,y)$, put
\begin{equation}\label{eq:highest}
 \istar=\max\{i:(\chi_i(e),x_i)\ne(\chi_i(f),y_i)\},
 \qquad u=q_{\istar}.
\end{equation}
Thus $\istar$ is an absolute word coordinate, whereas $u$ is the size of the suffix strictly below it. For example, take radices $(n_0,n_1,n_2)=(3,4,5)$ and fix a word-block label. Writing the letters from highest to lowest, the words $x=(2,1,3)$ and $y=(2,4,1)$ have $\istar=1$ and $u=q_1=3$; their difference at position $0$ belongs to the lower-digit tail. If instead $y=(4,1,3)$, then $\istar=2$ and $u=q_2=12$. A difference of seed color also counts as a difference at its position, even when the numerical letters agree. The separation estimates below show when the highest discrepancy dominates the tail.
Let $R_0=\dot\bigcup_{e\in\Gamma}W_e$.  For $0\le a\le P$ put
\begin{equation}\label{eq:proj}
 \Pi_a=I_{Q/q_a}\otimes J_{q_a};\qquad \Pi_0=I_Q,\quad \Pi_P=J_Q.
\end{equation}
The zero--one matrix $\Pi_a$ is the \emph{connection kernel at level $a$}; its row sum is $q_a$. The kernels satisfy
\begin{equation}\label{eq:projrule}
 \Pi_a\Pi_b=q_{\min(a,b)}\Pi_{\max(a,b)}.
\end{equation}
Let $T_s$ be the permutation matrix with $(T_s)_{ef}=1$ when $f=e+s$ and $0$ otherwise. Choose a set $\mathcal S\subset\Gamma\setminus\{0\}$ of pairwise distinct vectors and assign each $s\in\mathcal S$ a level $a_s\in\{0,\ldots,P\}$.  Define
\begin{equation}\label{eq:R}
 A_R=\bigoplus_{e\in\Gamma}A_{W_e}+\sum_{s\in\mathcal S}T_s\otimes \Pi_{a_s}.
\end{equation}
Because no nonzero displacement is repeated, $R$ is simple.  Its degree and order are exactly
\begin{equation}\label{eq:degree}
 \omega=D+\sum_{s\in\mathcal S}q_{a_s},\qquad N=2^\ell Q,
 \qquad r=N-1-\omega.
\end{equation}
Thus $\overline R$ is $r$-regular. A level-$a$ generator contributes exactly $q_a$ new neighbors to each vertex. Choosing how many generators to use at each level therefore encodes the desired ordinary degree.

For example, in a two-position word the suffix sizes are $q_0=1$, $q_1=n_0$ and $q_2=n_0n_1$. Between two blocks, a level-$0$ connection is a matching. A level-$1$ connection keeps the upper letter fixed and allows every lower letter, adding $n_0$ neighbors per vertex. A level-$2$ connection joins the two whole blocks and adds $Q$ neighbors per vertex. The same rule defines the connections at any depth.

For the main construction, the lower part has $p$ positions with seed $H_0$, and the upper part has $k$ positions using the $A/B$ alphabet. We call the lower part the pad. Set
\begin{equation}\label{eq:mixed}
 b=31250=2(125)^2,\quad n=125,\quad P=p+k,\quad
 Q=b^pn^k,\quad N=2^\ell Q.
\end{equation}
The absolute position $i$ and upper-position index $j$ are related by $i=p+j$. The three regions of word positions are
\[
\underbrace{0\le i<p}_{\text{pad }H_0},\qquad
\underbrace{i=p+j,\ 0\le j<p}_{\text{unamplified upper prefix}},\qquad
\underbrace{i=p+j,\ p\le j<k}_{\text{amplified upper positions}}.
\]
The boundary $j=p$ leaves a prefix of scale $b^pn^p$ without amplification; this scale will be negligible in the high-girth comparison, while complete-subspace amplification is reserved for higher positions.
At the $j$-th upper position, choose $A$ or $B$ according to the label bit $e_j$ if $j<\ell$, and use $A$ if $j\ge\ell$. All pad positions use $H_0$. The label-dependent color is part of the letter, so the highest unequal letter distinguishes any two roots when $\ell<k$.

Accordingly, $\chi_i=0$ for $i<p$, $\chi_{p+j}(e)=e_j$ for $j<\ell$, and $\chi_{p+j}=0$ for $j\ge\ell$. We retain $\istar$ from Eq.~\eqref{eq:highest} for the highest unequal absolute position. Unless otherwise specified, implicit constants may depend on the fixed pattern, seed system and baseline size $M=M(F)$ chosen before $r$, but not on $r,p,k,\ell$ or the compared roots.

At low proper levels we use only vectors supported in label bits $1,\ldots,p-1$. For the amplified upper levels $a=p+j$ with $j\ge p$, we will reserve $M$ vectors forming a full binary subspace with its zero element deleted.  Remaining level slots encode a radix digit.  Top-level generators use the full-block kernel $\Pi_P$ and are chosen from the disjoint coset with bit $0=1$; their number need not be bounded independently of $r$.  Lemma~\ref{lem:subspaces} and the capacity estimate in Section~\ref{sec:arithmetic} give explicit conditions for choosing these vectors without repetition.

\paragraph{Principal scale notation}
The size variables and bounds used below have the following meanings:
\begin{center}
\small
\begin{tabular}{@{}ll@{}}
\toprule
Symbol & Meaning \\
\midrule
$Q$, $N=2^\ell Q$ & order of a word block; total host order \\
$D$, $\omega$, $r$ & word-block degree; degree of $R$; degree of $\overline R$ \\
$\istar$, $u=q_{\istar}$ & highest differing word coordinate; lower-suffix size \\
$\varepsilon=\omega/N$ & sparsity parameter in the quotient bounds \\
\bottomrule
\end{tabular}
\end{center}
\Needspace{14\baselineskip}
\section{Shell invariance and local count separation}\label{sec:digits}
In the lexicographic words, rooted statistics form positional expansions. For two different roots, the highest unequal coordinate must dominate both the lower digits and the contributions of added connections. Equal suffix degrees across seed colors make this possible: first for triangles, and then for four-cycles in the bipartite construction.
\subsection{A common connection identity}
The common structural feature is a nested algebra of suffix relations: a connection may use several disjoint suffix shells without disturbing the differences between rooted triangle counts. The following statement includes both the binary labels used in the classification theorem and the cyclic labels used in the quantitative refinement. Retain the radices $n_i$, suffix sizes $q_i$ and kernels $\Pi_i$ from Section~\ref{sec:wiring}. Define
\begin{equation}\label{eq:shells}
 E_0=I_Q,\qquad E_i=\Pi_i-\Pi_{i-1}\quad(1\le i\le P).
\end{equation}
These symmetric zero--one matrices have disjoint supports partitioning $J_Q$, and their row sums are
\begin{equation}\label{eq:shell-degrees}
 \operatorname{rowSum}(E_0)=1,\qquad \operatorname{rowSum}(E_i)=q_i-q_{i-1}\quad(1\le i\le P).
\end{equation}
Equation~\eqref{eq:projrule} makes the span of the $\Pi_i$ a commutative algebra, all of whose elements have constant diagonal. Moreover, $\Pi_i/q_i$ is the averaging projection onto vectors constant on each suffix class of size $q_i$. A \emph{shell connection} is a sum of distinct $E_i$; disjointness makes it a zero--one matrix even when several levels are used.

\begin{lemma}[Uniform triangle increment]\label{lem:shell-invariance}
Let $\Lambda$ be a finite abelian group. For each $e\in\Lambda$, let $W_e$ be a lexicographic word graph with the same radices and the same seed degree at each position, but with arbitrary seed colors. Let $D_i$ be the common degree of a suffix of order $q_i$. Choose a set $\mathcal S\subseteq\Lambda\setminus\{0\}$ closed under negation, and assign to each $s\in\mathcal S$ a shell connection $K_s$, with $K_{-s}=K_s$. Then
\begin{equation}\label{eq:shell-host}
 A_R=\bigoplus_{e\in\Lambda}A_{W_e}+\sum_{s\in\mathcal S}T_s\otimes K_s
\end{equation}
defines a simple regular graph, and there is a constant $C$, independent of $e,x$ and the color assignment, such that
\begin{equation}\label{eq:shell-increment}
 \tau_R(e,x)=\tau_{W_e}(x)+C.
\end{equation}
Thus all differences between the rooted triangle counts of the word graphs are preserved.
\end{lemma}
\begin{proof}
Put $B_e=A_{W_e}$. The induced graph on each suffix class of order $q_i$ is $D_i$-regular, regardless of its prefix and color. Consequently
\begin{equation}\label{eq:shell-root}
 (B_e\Pi_i)_{xx}=D_i,\qquad
 (\Pi_iB_e\Pi_j)_{xx}=q_{\min(i,j)}D_{\max(i,j)}.
\end{equation}
For $i\le j$, the second expression counts the edges from the $q_i$ vertices in the smaller suffix class into the enclosing $q_j$-vertex class: each has $D_j$ neighbors there. The other case follows by symmetry.

The algebraic argument below applies to every finite abelian label group, including groups with elements of order two. It requires the suffix classes to have the same degrees across colors, but does not require the seed adjacency matrices themselves to commute.

A closed three-step walk explains the relevant connection types. One cross-block step cannot return to its starting block; two steps must have opposite displacements $s,-s$; three must have displacements summing to zero. The suffix algebra will show that both kinds of returning cross-block walks contribute the same amount at every root, independently of the seed colors.

Every shell connection belongs to the span of the $\Pi_i$; in particular its coefficients in $K=\sum_i b_i\Pi_i$ are well defined. Put
\begin{equation}\label{eq:shell-phi}
 \Phi(K)=\sum_{i,j}b_ib_jq_{\min(i,j)}D_{\max(i,j)}.
\end{equation}
Equations~\eqref{eq:projrule} and~\eqref{eq:shell-root} imply
\[
 (B_eK^2)_{xx}=(KB_fK)_{xx}=\Phi(K)
\]
for all $e,f,x$. For three such matrices $K_s,K_t,K_u$, the constant diagonal of their product is
\begin{equation}\label{eq:shell-psi}
 \Psi(K_s,K_t,K_u)
 =\sum_{i,j,h}b_{s,i}b_{t,j}b_{u,h}
     q_{\min(i,j,h)}q_{\operatorname{med}(i,j,h)},
\end{equation}
where $\operatorname{med}$ denotes the middle index.

Expand the diagonal of $A_R^3$ in ordered matrix words. The one-connection terms vanish by the block-level observation above. For two connections of displacements $s,-s$, the three placements of the within-block factor contribute $3\Phi(K_s)$ after summation over $s\in\mathcal S$; the count also covers $s=-s$. Three-connection terms are indexed by $s+t+u=0$. Thus
\begin{equation}\label{eq:shell-exact}
 \begin{split}
 2\bigl(\tau_R(e,x)-\tau_{W_e}(x)\bigr)
 &=3\sum_{s\in\mathcal S}\Phi(K_s)\\
 &\quad+\sum_{\substack{s,t,u\in\mathcal S\\s+t+u=0}}
          \Psi(K_s,K_t,K_u),
 \end{split}
\end{equation}
which is independent of the root and colors. Every nonzero displacement is used once; different displacements and the within-block edges have disjoint supports. Symmetry follows from $K_{-s}=K_s$, and the common row sum gives regularity.
\end{proof}

We apply Lemma~\ref{lem:shell-invariance} with binary label groups in the classification proof and with cyclic labels in the quantitative triangle construction. In both cases the shell connections adjust the common degree while preserving every rooted triangle-count difference.

\subsection{Triangle and four-cycle separation}
The bounded proper-level generator budget in the four-cycle statement below is uniform in word length; this is the form used by the localization estimate.
\begin{lemma}[Triangle invariance]\label{lem:trinvariance}
With $R$ as in~\eqref{eq:R}, all new terms in $\diag(A_R^3-A_{R_0}^3)$ are constant across the vertices.  In particular,
\begin{equation}\label{eq:trinvariance}
 \Delta\tau_R=\Delta\tau_{R_0}
\end{equation}
for every two roots.
\end{lemma}
\begin{proof}
Apply Lemma~\ref{lem:shell-invariance} with $\Lambda=\mathbb F_2^\ell$ and $K_s=\Pi_{a_s}=E_0+\cdots+E_{a_s}$. Every displacement is its own negative, and the original construction is exactly~\eqref{eq:shell-host} with these choices.
\end{proof}

For the words in~\eqref{eq:mixed}, Eq.~\eqref{eq:trianglelex} gives
\begin{equation}\label{eq:triword}
 \tau_{R_0}(e,x)=C+\sum_{i=0}^{P-1}q_i^2\tau_{S_{i,\chi_i(e)}}(x_i).
\end{equation}
Here $C$ depends only on the common degrees and radices. Compare two roots at their highest unequal coordinate $\istar$ and set $u=q_{\istar}$. The associated seed counts differ by a nonzero integer, contributing at least $u^2$ in absolute value. At a higher-alphabet coordinate, the lower upper digits have total span at most $932u^2/(125^2-1)$, and the lower padding digits contribute at most $14\,563\,432u^2/(31250^2-1)$. If the highest difference lies in the padding itself, only the latter geometric tail occurs. Thus, with
\begin{equation}\label{eq:eta}
 \eta=\frac{932}{125^2-1}+\frac{14\,563\,432}{31250^2-1}<0.075,
\end{equation}
we obtain
\begin{equation}\label{eq:trigap}
 |\Delta\tau_R|\ge (1-\eta) u^2.
\end{equation}
The joint distinctness of each seed alphabet ensures that the same argument applies to different block labels: a differing label bit is represented by a differing upper seed color. In particular, the bound is uniform over all pairs of roots.

\begin{lemma}[Four-cycle separation under layered connections]\label{lem:c4}
Take $k$ positions from the two $1250$-vertex seeds $U,V$ of Lemma~\ref{lem:seeds}, with $\ell\le k$, using color $e_i$ at position $i<\ell$ and fixed color $0$ thereafter.  Use any collection of distinct connection generators with a fixed bounded count at each proper level.  Then $\tau_R$ is constant, and for every pair of roots
\begin{equation}\label{eq:c4inv}
 \Delta c_4^R=\Delta c_4^{R_0},\qquad
 |\Delta c_4^R|\ge\left(1-\frac{834\,684}{1250^3-1}\right)u^3.
\end{equation}
\end{lemma}
\begin{proof}
The seeds $U,V$ are bipartite, so their triangle counts are zero. Equation~\eqref{eq:trianglelex} makes $\tau_{R_0}$ constant, and Lemma~\ref{lem:trinvariance} then makes $\tau_R$ constant. For any regular simple graph $X$,
\begin{equation}\label{eq:c4walk}
 (A_X^4)_{vv}=2c_4^X(v)+2d(X)^2-d(X).
\end{equation}
In the fourth-order expansion of Appendix~\ref{app:fourth}, the
potentially color-dependent alternating two-connection terms reduce to
mixed same-label common-neighbor counts. The all-depth formula in
Lemma~\ref{lem:mixed-cn} shows that these depend only on the block
displacement, not on the root: same-color seed positions contribute
$310$ and opposite-color positions contribute zero. The remaining
connection words also have root-independent diagonals, as verified
term by term in that appendix. Equation~\eqref{eq:c4walk} therefore gives the first identity in~\eqref{eq:c4inv}.

For the digit bound, expanding $A_{X[Y]}=A_X\otimes J_q+I\otimes A_Y$ gives
\[
 c_4^{X[Y]}(x,y)=q^3c_4^X(x)+4q^2d_Y\tau_X(x)+c_4^Y(y)+C,
\]
where $C$ depends only on orders and degrees. Each outer seed is $U$ or $V$, so the triangle term vanishes. Induction gives
\[
 c_4^{R_0}(e,x)=C_k+\sum_i1250^{3i}c_4^{S_{i,\chi_i(e)}}(x_i).
\]
The highest unequal letter is at $\istar$ and has gap at least one. The lower tail is at most $834\,684/(1250^3-1)$ times its weight $u^3$. This proves the second assertion.
\end{proof}

\section{Rooted localization by block centering}\label{sec:local}
For arbitrary connected rooted patterns, align the two roots and subtract a common blockwise average. The cancellation mechanism is a leaf sum: a nonempty forest of centered edges has a leaf other than the prescribed root, and, after fixing the other images, summing over that leaf gives a zero row or column sum. The remaining uncentered factors are constant on the relevant blocks. Hence only markings containing a cycle can distinguish two roots. The localization bound below makes this observation uniform even when the number of full-block connections grows.

\begin{lemma}[Common block average after root alignment]\label{lem:alignment}
Let $R$ be an auxiliary graph constructed in Section~\ref{sec:wiring}. Fix roots $v=(e,x)$ and $w=(f,y)$, and let $\istar$ be their highest unequal absolute coordinate. Partition $V(R)$ by the word-block label and the letters strictly above $\istar$; every block has size $L_*=q_{\istar+1}$. There is a vertex permutation $\varphi$ with $\varphi(v)=w$ such that the relabeled adjacency matrix
\[
 A'_R(z,z')=A_R(\varphi(z),\varphi(z'))
\]
has the same blockwise average $\mathcal U$ as $A_R$. Write $K=A_R-\mathcal U$ and $K'=A'_R-\mathcal U$, and put
\[
 \omega_*=D_{\istar+1}+\sum_{s:a_s\le\istar}q_{a_s}.
\]
Then every block of $K$ and $K'$ has zero row and column sums,
\[
 0\le\mathcal U(z,z')\le1,\quad |K(z,z')|,|K'(z,z')|\le1,
 \quad \sum_{z'}\mathcal U(z,z')=\omega,
\]
and
\[
 \max\{\|K\|_{\infty\to\infty},\|K'\|_{\infty\to\infty}\}
 \le2\omega_*,\qquad \omega_*\le\omega.
\]
There is no restriction on the number of top-level connections.
\end{lemma}
\begin{proof}
Put $s=e+f$. Equality of the two colors at positions $i>\istar$ gives $\chi_i(s)=0$. Translate every word-block label by $s$ and, at each position $i\le\istar$, apply a fixed permutation of $[n_i]$ taking $x_i$ to $y_i$. Apply the same alphabet permutation in every word block. The resulting permutation $\varphi$ sends $v$ to $w$ and preserves all equality relations that define the connection kernels.

For blocks $\mathfrak b,\mathfrak c$, let $c_{\mathfrak b\mathfrak c}$ be the common row sum of $A_R[\mathfrak b,\mathfrak c]$. These blocks are biregular. A lexicographic edge from a position above $\istar$ contributes either $0$ or $J_{L_*}$, and its value is unchanged by $\varphi$. Within a word block and a fixed upper prefix, lower edges form a regular suffix graph of degree $D_{\istar+1}$. A connection of level $a>\istar$ contributes either $0$ or $J_{L_*}$ between the relevant blocks; one of level $a\le\istar$ contributes row and column sum $q_a$. Translation of the word-block labels preserves each displacement. Thus these block row sums agree entrywise for $A_R$ and $A'_R$.

Define
\[
 \mathcal U[\mathfrak b,\mathfrak c]
 =\frac{c_{\mathfrak b\mathfrak c}}{L_*}J_{L_*}.
\]
This proves the common-average assertion and the blockwise zero sums. The entry bounds follow because $A_R,A'_R$ are zero--one matrices and $\mathcal U$ is their block average. Its row sum is $\omega$. All all-ones blocks disappear upon centering, so only suffix edges and connections of level at most $\istar$ contribute to $K$ or $K'$. Centering a nonnegative block of row sum $c$ gives absolute row sum at most $2c$. Summing these contributions proves the asserted bound.
\end{proof}

\Needspace{17\baselineskip}
\begin{theorem}[Rooted localization]\label{thm:local}
Fix regular seed alphabets with orders $2\le n_i\le n_{\max}$ and degrees independent of color at each position. Suppose each proper connection level $a<P$ carries at most $M_*$ generators; the number at level $P$ is unrestricted. Let $v,w$ be roots of an auxiliary graph $R$ from Section~\ref{sec:wiring}, with highest unequal coordinate $\istar$, and set $u=q_{\istar}$. For any fixed connected rooted graph $(H,a)$ of order $h_H$:
\begin{enumerate}[label=(\roman*)]
\item if $H$ is a tree, $\Hom_v(H,a;R)=\Hom_w(H,a;R)$;
\item if $H$ contains a cycle, then
\begin{equation}\label{eq:loc}
 |\Delta\Hom(H,a;R)|
 \le C_{H,M_*}\omega^{h_H-g(H)}u^{g(H)-1}.
\end{equation}
\end{enumerate}
The constant may depend on the fixed alphabets, $H$ and $M_*$, but not on $P,\ell$, the selected generators or the roots.
\end{theorem}
\begin{proof}
Use Lemma~\ref{lem:alignment} to compare $\mathcal U+K$ and $\mathcal U+K'$ at $v$. Since the radices are bounded and at least two, geometric summation gives
\begin{equation}\label{eq:Knorm}
 \omega_*=D_{\istar+1}+\sum_{a_s\le\istar}q_{a_s}
 \le C_{M_*}u,\qquad \omega_*\le\omega.
\end{equation}
Top-level connections occur only in $\mathcal U$.

\emph{Expansion and cancellation.}
Expand each source edge weight as $\mathcal U+K$ and fix the averaging block assigned to each source vertex. Let $S$ be the edges marked $K$. If $S$ is a nonempty forest, choose any component containing a marked edge. That component is a nontrivial tree and has at least two leaves; at most one is the distinguished source root. Hence some marked leaf $z$ is not the root, even when the marked forest has several components or isolated vertices. Fix the averaging blocks $\mathfrak b_x$ of all source vertices and the images in $R$ of every vertex except $z$. Write $z'$ for the unique neighbor of $z$ along an edge of $S$. All other edges incident to $z$ are marked $\mathcal U$, and their weights depend only on the fixed block labels. Thus the sum of the corresponding term over the image $y$ of $z$ contains the factor
\[
 \left(\prod_{zz''\in E(H)\setminus S}
       \frac{c_{\mathfrak b_z\mathfrak b_{z''}}}{L_*}\right)
 \sum_{y\in\mathfrak b_z}K(y,\phi(z'))=0,
\]
where $\phi(z')$ is already fixed. The last equality is the zero column sum of the indicated block of $K$; reversing the orientation uses its zero row sum. Hence every nonempty forest marking vanishes, for $K$ and for $K'$. The empty marking depends only on the common $\mathcal U$ and cancels in the root difference. Consequently every surviving marking contains a cycle of $H$.

\emph{The spanning-tree bound.}
Choose $g(H)-1$ consecutive edges on such a cycle and extend them to a spanning tree of $H$. The selected edges are marked $K$ (or $K'$). Assign the root image and sum the other images along the tree. Each selected edge has absolute row sum at most $C_{M_*}u$; every other tree edge has absolute row sum at most $2\omega$. Non-tree factors have absolute value at most one. A marking is therefore bounded by
\[
 (C_{M_*}u)^{g(H)-1}(2\omega)^{h_H-g(H)}.
\]
There are finitely many markings, depending only on $H$, so their sum gives Eq.~\eqref{eq:loc}. A tree has only empty or forest markings; equivalently, its rooted homomorphism count is $\omega^{h_H-1}$.
\end{proof}

For the concrete alphabets used here, $d_{\istar}\ge1$ and hence
\begin{equation}\label{eq:uomega}
 u\le d_{\istar}q_{\istar}\le D\le\omega.
\end{equation}
For any integer $q$ with $3\le q\le g(H)$, Eq.~\eqref{eq:loc} consequently implies
\begin{equation}\label{eq:locq}
 |\Delta\Hom(H,a;R)|\le C_{H,M_*}\omega^{h_H-q}u^{q-1}.
\end{equation}
Indeed, the ratio of the first bound to the second is $(u/\omega)^{g(H)-q}\le1$.

\section{Complement expansions and leading cyclic quotients}\label{sec:quotients}
Expanding an $F$-degree in $\overline R$ yields injective counts of edge supports of $F$. Partition inversion expresses these in rooted homomorphism counts. The localization bound then restricts which cyclic quotient types can contribute at leading order.

Throughout this section, $R$ is an auxiliary graph constructed in Section~\ref{sec:wiring} and subject to the hypotheses of Theorem~\ref{thm:local}; its seed degrees are positive. Fix roots $v,w$ and write $u=q_{\istar}$; we use the sparsity parameter $\varepsilon$ from Section~\ref{sec:wiring}. For a nonempty edge set $E'\subseteq E(F)$, let $J=F[E']$ be its edge support after removing all isolated vertices, and put $s=|V(J)|$.  Define
\[
 I_J(v)=\sum_{a\in V(J)}\Inj_v(J,a;R).
\]
\begin{lemma}[Rooted complement identity]\label{lem:complement}
If $G=\overline R$ has order $N\ge2h$, where $h=|V(F)|$, then for any two roots in $G$,
\begin{equation}\label{eq:complement}
 a_F\Delta d_F^G
 =\sum_{\varnothing\ne E'\subseteq E(F)}(-1)^{|E'|}
 (N-s-1)_{h-s}\,\Delta I_J,
\end{equation}
where $(x)_m=x(x-1)\cdots(x-m+1)$ for $m\ge1$ and $(x)_0=1$.
\end{lemma}
\begin{proof}
Expand each adjacency indicator of $G$ in the injective count~\eqref{eq:rootedinj} as $1-A_R$.  If $s=h$, the distinguished source root is necessarily in the support and the coefficient is $1=(N-s-1)_0$.  Hence assume $s<h$ in the following calculation.  For fixed nonempty $E'$ the vertices outside its support carry no edge restrictions.  When the distinguished source root belongs to $J$, these other $h-s$ vertices can be injected in $(N-s)_{h-s}$ ways.  When the root lies outside $J$, subtract the contribution of embeddings of $J$ whose images meet the distinguished vertex of $R$.  This second contribution has multiplicity $(h-s)(N-s-1)_{h-s-1}$.  The coefficient is therefore
\[
 (N-s)_{h-s}-(h-s)(N-s-1)_{h-s-1}=(N-s-1)_{h-s}.
\]
The contribution of $E'=\varnothing$ is independent of the root.  This proves~\eqref{eq:complement}.
\end{proof}

For a graph $J$, injective rooted homomorphisms can be expressed by partition-lattice inversion:
\begin{equation}\label{eq:mobius}
 \Inj_v(J,a;R)=\sum_{\pi\in\Pi_{\rm ind}(J)}\mu(\pi)
 \Hom_v(J/\pi,[a];R).
\end{equation}
Here $\Pi_{\rm ind}(J)$ consists of partitions whose blocks are independent in $J$, with Möbius coefficient
\[
 \mu(\pi)=\prod_{B\in\pi}(-1)^{|B|-1}(|B|-1)!.
\]
Parallel quotient edges may be merged because the target adjacency is zero--one. In summing over source roots, a quotient vertex corresponding to a block $B$ retains its multiplicity $|B|$.

\begin{lemma}[Quotient-order bound]\label{lem:quotient}
Fix $F$ and two roots at suffix scale $u$. Let a quotient of an edge support $J$ have $m$ vertices, $c$ nontrivial connected components and identification defect $d=s-m$. Let $H$ be the quotient component containing the distinguished root. If $H$ has a cycle, choose an integer $q$ with $3\le q\le g(H)$. The contribution of this rooted quotient to~\eqref{eq:complement} is bounded by
\begin{equation}\label{eq:quotient}
 C_FN^{h-q-d}\left(\frac\omega N\right)^{m-c-q+1}u^{q-1},
 \qquad m-c-q+1\ge0.
\end{equation}
A quotient with the root in a tree component makes no contribution to the root difference.
\end{lemma}
\begin{proof}
Write the quotient as a disjoint union $J/\pi=H_1\dot\cup H_2\dot\cup\cdots\dot\cup H_c$, with the distinguished quotient root $[a]$ in $H_1$. For each root $v\in V(R)$, homomorphisms into $R$ factor exactly:
\[
 \Hom_v(J/\pi,[a];R)=\Hom_v(H_1,[a];R)
       \prod_{j=2}^{c}\Hom(H_j;R).
\]
The factors outside the rooted component are independent of $v$. Therefore
\[
 \Delta\Hom(J/\pi,[a];R)=\Delta\Hom(H_1,[a];R)
       \prod_{j=2}^{c}\Hom(H_j;R).
\]
If $H_1$ is a tree, then $\Hom_v(H_1,[a];R)=\omega^{|V(H_1)|-1}$ is independent of $v$, and the difference is zero. If $H_1$ contains a cycle, use~\eqref{eq:locq} on $H_1$ and the spanning-tree bound $\Hom(H_j;R)\le N\omega^{|V(H_j)|-1}$ on each other component. This gives
\[
\begin{aligned}
 |\Delta\Hom(J/\pi,[a];R)|
 &\le C_F\omega^{|V(H_1)|-q}u^{q-1}
 \prod_{j=2}^{c}N\omega^{|V(H_j)|-1},\\
 &= C_FN^{c-1}\omega^{m-c-q+1}u^{q-1}.
\end{aligned}
\]
Multiplication by the coefficient $O_F(N^{h-s})$ from~\eqref{eq:complement} proves~\eqref{eq:quotient}. The rooted cyclic component has at least $q$ vertices, and every other nontrivial component has at least two, so $m-c-q+1\ge0$.
\end{proof}

We distinguish the \emph{identification defect} $d=s-m$ from the \emph{sparsity exponent} $m-c-q+1$ in~\eqref{eq:quotient}. Vertex identification may create a triangle even when $F$ is triangle-free. The localization bound therefore uses the girth of the rooted quotient component.

A \emph{discrete quotient} has only singleton partition blocks, although its edge support need not be connected. Which types contribute at leading order depends on the structure of $F$. A type is called \emph{zero-loss} at a specified leading scale when its quotient bound retains that scale without an additional factor of $N^{-1}$ or $\varepsilon$. For a concrete example, write $P_4=a-b-c-d$: identifying the nonadjacent endpoints $a,d$ gives the quotient triangle $[a=d]-b-c-[a=d]$. Thus a non-star tree can produce a rooted triangle-count term after partition inversion, even though the original pattern has no cycle. Table~\ref{tab:leading} displays the contributions; Lemma~\ref{lem:leading} accounts for the precise source-root multiplicities.
\begin{table}[htbp]
\centering
\caption{Leading contributions to $a_F\Delta d_F^{\overline R}$ under the indicated pattern hypotheses. The four-cycle row requires triangle-free $F$ and constant $\tau_R$.}
\label{tab:leading}
\small
\begin{tabular}{@{}lll@{}}
\toprule
Support and identification & Main scale & Contribution\\
\midrule
Triangle, no identification & $N^{h-3}u^2$ & $-6c_3(F)N^{h-3}\Delta\tau_R$\\
$P_4$, endpoints identified & $N^{h-4}u^2$ & $8p_4(F)N^{h-4}\Delta\tau_R$\\
Four-cycle, no identification & $N^{h-4}u^3$ & $8c_4(F)N^{h-4}\Delta c_4^R$\\
Shortest $g$-cycle, $g\ge5$ & $N^{h-g}u^{g-1}$ & $(-1)^ggc_g(F)N^{h-g}\Delta z_g^R$\\
\bottomrule
\end{tabular}
\end{table}

The loss mechanism is controlled by two quantities: the identification defect $d=s-m$ and the exponent $A_q=m-c-q+1$ from Lemma~\ref{lem:quotient}. A rooted tree quotient cancels exactly. For a rooted cyclic quotient, $A_q=0$ forces a single chordless cycle; an identification lowers the power of $N$ relative to terms of the same cyclic girth, but the endpoint-identified $P_4$ still supplies a separate leading scale. The only triangle-producing identification in a triangle-free four-vertex support closes the endpoints of $P_4$. Lemma~\ref{lem:leading} proves the resulting classification; Appendix~\ref{app:quotient-ledger} records the source-root multiplicities and loss factors.

For triangle-free patterns of girth $g\ge5$, the contracted $P_4$ and the uncontracted shortest cycle must be retained at \emph{different} principal scales. Neither is automatically an error relative to the other. Their possible cancellation is the reason for the spectral argument in Section~\ref{sec:spectral}.

\begin{lemma}[Leading types]\label{lem:leading}
Let $F$ be connected and $h=|V(F)|$.
\begin{enumerate}[label=(\alph*)]
\item If $F$ contains a triangle, the only root-dependent zero-loss type at scale $N^{h-3}u^2$ is an uncontracted triangle support. Summing over these supports gives coefficient $-6c_3(F)$ multiplying $\Delta\tau_R$.
\item If $F$ is triangle-free and contains a $P_4$, the only root-dependent zero-loss type among non-discrete quotients at scale $N^{h-4}u^2$ identifies the endpoints of a $P_4$. Its summed coefficient is $8p_4(F)$ multiplying $\Delta\tau_R$.
\item If $F$ has finite girth $g\ge5$, the only zero-loss type among discrete cyclic quotients at scale $N^{h-g}u^{g-1}$ is a single $C_g$ support. Their sum has coefficient $(-1)^ggc_g(F)$ multiplying $\Delta z_g$.
\end{enumerate}
In (a), all other root-dependent types have an additional factor $\omega/N$ or $N^{-1}$ relative to the triangle scale. The same loss holds in (b) for the other non-discrete types relative to the path scale, and in (c) for the other discrete cyclic types relative to the shortest-cycle scale.
\end{lemma}
\begin{proof}
Write $\mathcal T$ for a fixed root-dependent quotient term and let $H$ be its rooted cyclic component. In this proof choose $q=g(H)$, the actual girth of that component, and put $A_q=m-c-q+1$.  By Lemma~\ref{lem:quotient}, a rooted cyclic
quotient contributes at most
\[
 C_FN^{h-q-d}\eps^{A_q}u^{q-1},
 \qquad A_q\ge0,\quad q\ge3,\quad u/N\le\eps.
\]
The root component has at least $q$ vertices. Notice that every other component of the edge support is nontrivial, hence has at least two vertices; a tree attached to the rooted cycle adds further vertices to its component. If $A_q=0$, then
$m=q+c-1$; since every other nontrivial component has at least two
vertices, this forces $c=1$ and $m=q$.  Since the edge support and all
its quotients have no isolated vertices, the quotient then consists
of exactly one connected cyclic component on $q$ vertices. At $q=g(H)$, that component must be the $q$-cycle itself: an additional chord would create a shorter cycle, contrary to the choice $q=g(H)$.

Relative to the triangle scale $N^{h-3}u^2$, the bound is
\begin{equation}\label{eq:triangle-loss}
 \frac{|\mathcal T|}{N^{h-3}u^2}
 \le C_FN^{-d}(u/N)^{q-3}\eps^{A_q}
 \le C_FN^{-d}\eps^{q-3+A_q}.
\end{equation}
The only zero-loss choice is $d=0$, $q=3$, $A_q=0$,
namely the uncontracted triangle.  It contributes $-6c_3(F)$, since
there are three root positions and $\Hom_v(K_3,a;R)=2\tau_R(v)$.

For a triangle-free source, compare instead with the contracted-path
scale $N^{h-4}u^2$:
\begin{equation}\label{eq:path-loss}
 \frac{|\mathcal T|}{N^{h-4}u^2}
 \le C_FN^{1-d}\eps^{q-3+A_q}.
\end{equation}
If $d=0$, the support cannot have $q=3$; when $q\ge4$, the
right-hand side in~\eqref{eq:path-loss} is not itself sufficient
for suppression because it contains a factor $N$.  For a discrete
support with $q\ge4$, we therefore retain its own scale
$N^{h-q}u^{q-1}$ rather than treating it as a contracted-path error.
Among non-discrete quotients ($d\ge1$), only
$d=1$, $q=3$, $A_q=0$ can have zero loss.
It has $s=m+d=4$, so precisely one independent pair of vertices
is identified in an edge support with four vertices and no isolates.
Such a triangle-free support is one of $2K_2$, $K_{1,3}$, $P_4$,
and $C_4$. An admissible identification in $2K_2$ or $K_{1,3}$
yields the tree $P_3$; in $C_4$ the identified vertices must be
opposite, again yielding $P_3$. In $P_4$, identification of a
nonadjacent pair at distance two also yields $P_3$, while identifying
the two endpoints yields $K_3$. Rooted homomorphism counts of $P_3$
are equal to $\omega^2$ at every vertex of an $\omega$-regular graph. Thus only
the endpoint identification in $P_4$ has a possibly nonconstant rooted
quotient count, namely $\Hom_v(K_3,a;R)=2\tau_R(v)$ for any choice of the triangle root $a$.
Its M\"obius sign $-1$ cancels the sign $(-1)^3$ of the three support
edges.  The four root positions contribute multiplicities
$2,1,1$, and $\Hom_v(K_3,a;R)=2\tau_R(v)$; the total coefficient is
therefore $8p_4(F)$.

If the source has girth $g\ge5$, any \emph{discrete} cyclic support
has rooted-component girth $q\ge g$.  Relative to the shortest-cycle
scale $N^{h-g}u^{g-1}$ its contribution is at most
\begin{equation}\label{eq:cycle-loss}
 \frac{|\mathcal T|}{N^{h-g}u^{g-1}}
 \le C_F(u/N)^{q-g}\eps^{A_q}
 \le C_F\eps^{q-g+A_q}.
\end{equation}
The only zero-loss choice is $q=g$, $A_q=0$;
by the preceding characterization, the support is a single chordless
$C_g$.  Its $g$ root positions give the coefficient
$(-1)^g g c_g(F)$ multiplying $\Delta z_g$.
Every other non-discrete term has at least one factor
$\eps$ or $N^{-1}$ relative to the path scale, and every other
discrete cyclic term has an extra factor $\eps$ relative to the
shortest-cycle scale.  If $F$ is a tree there are no discrete cyclic
supports.  The finite falling-factorial correction in
\eqref{eq:complement} is $1+O_F(N^{-1})$ relative to $N^{h-s}$,
which contributes only the stated lower-order terms.
\end{proof}

\paragraph{Uniform remainder bounds.}
The remainder is a finite sum over the edge supports, admissible partitions and source-root positions determined by the fixed graph $F$. In the triangle-free case, a non-discrete
cyclic quotient with $d\ge2$ is $O_F(N^{-1})$ relative to
$S_P=N^{h-4}u^2$. A non-discrete quotient with $d=1$ other than the
endpoint-identified $P_4$ satisfies $q-3+A_q\ge1$ and is
$O_F(\varepsilon)$ relative to $S_P$. For discrete cyclic quotients of a
source of girth $g\ge5$, those other than an uncontracted shortest cycle
satisfy $q-g+A_q\ge1$ and are $O_F(\varepsilon)$ relative to
$S_C=N^{h-g}u^{g-1}$. Here we use only $u/N\le\varepsilon$ and
$N^{-1}\le1$. Rooted-tree terms vanish identically, while the
falling-factorial coefficients introduce $O_F(N^{-1})$ relative errors.
Consequently a single constant $C_F$, independent of the prescribed
regularity and the two roots, bounds the combined remainder in terms
of $S_P+S_C$.

As a consequence, for triangle-free $F$ of girth $g\ge5$,
\begin{equation}\label{eq:principal}
\begin{split}
 a_F\Delta d_F^{\overline R}={}&8p_4(F)N^{h-4}\Delta\tau_R
 +(-1)^g gc_g(F)N^{h-g}\Delta z_g^R+E_F,\\
 |E_F|\le{}&C_F(\eps+N^{-1})
 \bigl(N^{h-4}u^2+N^{h-g}u^{g-1}\bigr).
\end{split}
\end{equation}
If $F$ contains triangles, the leading term is $-6c_3(F)N^{h-3}\Delta\tau_R$ with an $O_F(\eps+N^{-1})$ relative error at scale $N^{h-3}u^2$. For a non-star tree the corresponding term is $8p_4(F)N^{h-4}\Delta\tau_R$. When $F=P_4$, this coefficient also follows from the exact identity in Lemma~\ref{lem:basic-identities}. For the four-cycle principal term, constancy of $\tau_R$ is needed because vertex identification may create triangles in a quotient.
\begin{lemma}[Four-cycle principal term]\label{lem:c4main}
Let $F$ be a connected triangle-free graph of order $h$ containing a four-cycle. Let $R$ be an $\omega$-regular auxiliary graph constructed as in Section~\ref{sec:wiring}, with positive seed degrees and with bounded proper-level generator counts as in Theorem~\ref{thm:local}. Suppose $\tau_R$ is constant. For distinct roots $v,w$, put $u=q_{\istar}$ and $\eps=\omega/N$. Then
\begin{equation}\label{eq:c4main}
 a_F\Delta d_F^{\overline R}
 =8c_4(F)N^{h-4}\Delta c_4^R
 +O_F\bigl((\eps+N^{-1})N^{h-4}u^3\bigr).
\end{equation}
\end{lemma}
\begin{proof}
For discrete quotients use~\eqref{eq:quotient} with $q=4$; the only zero-loss cyclic support is a single $C_4$.  Its rooted homomorphism count is $z_4=2c_4+2\omega^2-\omega$, so summing the four root positions yields coefficient $8$.

For non-discrete quotients use the weaker $q=3$ bound and let $d\ge1$ be the identification defect.  Relative to the scale $N^{h-4}u^3$, the bound is
\begin{equation}\label{eq:c4ratio}
 C_FN^{1-d}\eps^{m-c-2}/u.
\end{equation}
If $d\ge2$ this is $O(N^{-1})$; if $d=1$ and $m-c-2\ge1$ it is $O(\eps)$.  In the remaining case $d=1$, $m=3,c=1,s=4$.  Since $F$ is triangle-free, the possible four-vertex supports are $P_4,K_{1,3},2K_2$, and $C_4$.  Identifying independent vertices gives a quotient of type $P_3$ or $K_3$.  Its rooted homomorphism counts in a regular graph are respectively $\omega^2$ and $2\tau_R$, both constant under the hypothesis.  The remaining terms have the advertised loss.
\end{proof}

\section{Spectral separation and amplification at large girth}\label{sec:spectral}
For a pattern of girth $g\ge5$, the complement expansion leaves two rooted statistics whose contributions may cancel: the triangle count $\tau_R$ and the shortest-cycle closed-walk count $z_g^R$. We isolate their differences at the highest unequal word coordinate. The starting observation is elementary: in a regular seed, the diagonal adjacency moments of orders zero, one and two are root-independent; the third moment, $2\tau$, is the first that can distinguish letters.

The top-level degree correction is invisible to differences of rooted spectral moments. Its word-coordinate factor is $J_Q$, which annihilates every zero-sum word layer. On the remaining constant-word layer, the induced Cayley operator has the same diagonal powers at every root. Thus these connections can adjust the degree without disturbing any $\Delta z_t^R$; Theorem~\ref{thm:spectral} makes this separation exact.

At an amplified coordinate, a complete binary subspace forms a Cayley clique on $M+1$ label vertices. In the $g$th moment, the term with three seed factors and $g-3$ label-level factors therefore carries the rooted triangle difference with a coefficient of order $M^{g-4}$. Terms involving higher seed moments lose at least one power of $M$. Three possible disturbances remain: higher moments at the distinguished seed, the geometric tail of earlier amplified coordinates, and the unamplified prefix. The proof controls these through $E_g$, an amplified-digit envelope $W_g$, and a prefix bound $P_g$, respectively. We first prove the exact spectral decomposition and then establish a positive margin. Throughout, $g$ is fixed before the amplification size and padding depth, and all are fixed before the degree or roots are chosen.

\begin{theorem}[Spectral digit decomposition]\label{thm:spectral}
Let $R$ be an auxiliary graph constructed in Section~\ref{sec:wiring}, and let $D_i$ be the suffix degree in Eq.~\eqref{eq:worddegree}. Define the matrix $\mathcal B_i$ on $\Gamma\times[n_i]$ by
\begin{equation}\label{eq:Bi}
 \mathcal B_i=D_iI+q_i\diag_{e\in\Gamma}(A_{S_{i,\chi_i(e)}})
 +\sum_{s:a_s\le i}q_{a_s}(T_s\otimes I_{n_i}).
\end{equation}
For every positive integer $t$ and every two roots $(e,x),(f,y)$ of $R$,
\begin{equation}\label{eq:spectral}
 \Delta z_t^R
 =\sum_{i=0}^{P-1}\frac1{q_i}
 \bigl[(\mathcal B_i^t)_{(e,x_i),(e,x_i)}-(\mathcal B_i^t)_{(f,y_i),(f,y_i)}\bigr].
\end{equation}
In particular, top-level generators do not affect these differences.
\end{theorem}
\begin{proof}
For each position write $\RR^{n_i}=\langle\one_i\rangle\oplus V_i^0$, where $V_i^0$ is the zero-sum subspace.  The word-vector space decomposes orthogonally into the all-constant subspace and the subspaces
\begin{equation}\label{eq:orth}
 \left(\bigotimes_{a>i}\RR^{n_a}\right)\otimes V_i^0\otimes
 \left(\bigotimes_{a<i}\langle\one_a\rangle\right).
\end{equation}
On the $i$-th such subspace, lower lexicographic edges act by $D_iI$, the edge at position $i$ acts by $q_iA_{S_{i,\chi_i(e)}}$, and higher edges act by zero because they contain an all-ones factor on $V_i^0$.  Similarly, $\Pi_{a_s}$ acts by $q_{a_s}I$ if $a_s\le i$ and by zero otherwise.  The resulting matrix is the restriction of $\mathcal B_i$ to $V_i^0$.

The diagonal of the orthogonal projection onto the constant lower letters is $1/q_i$. To recover the restriction to $V_i^0$ from $\mathcal B_i$, subtract the one-dimensional constant-letter part. Since every $S_{i,\chi_i(e)}$ has degree $d_i$, the induced operator on $\Gamma$ is
\[
 \mathcal C_i=(D_i+q_id_i)I_\Gamma
       +\sum_{s:a_s\le i}q_{a_s}T_s.
\]
Its contribution to the rooted diagonal of the $i$th layer is
$-\,(\mathcal C_i^t)_{ee}/(q_in_i)$. Each $\mathcal C_i$ is a Cayley convolution: it commutes with every $T_h$, so $(\mathcal C_i^t)_{ee}$ is independent of the word-block label $e$. The entirely constant word subspace similarly contributes the diagonal of
\[
 \mathcal C_P=DI_\Gamma+\sum_{s\in\mathcal S}q_{a_s}T_s
\]
with coefficient $1/Q$; its diagonal is also independent of $e$. In fact, the complete diagonal identity is
\begin{equation}\label{eq:fullspectral}
 z_t^R(e,x)=\frac{(\mathcal C_P^t)_{ee}}{Q}
 +\sum_{i=0}^{P-1}\left[
 \frac{(\mathcal B_i^t)_{(e,x_i),(e,x_i)}}{q_i}
 -\frac{(\mathcal C_i^t)_{ee}}{q_in_i}\right].
\end{equation}
Both families of Cayley corrections are independent of the root's word-block label; subtracting the formula for two roots gives~\eqref{eq:spectral}. Top-level $J_Q$ annihilates every zero-sum layer~\eqref{eq:orth}; its only contribution appears through $\mathcal C_P$ and thus cancels.
\end{proof}

\begin{lemma}[Dependence on the local color]\label{lem:localcolor}
For each absolute position $i$ and each $t\ge0$, the value
$(\mathcal B_i^t)_{(e,x_i),(e,x_i)}$ depends on $e$ only through $\chi_i(e)$.
Consequently, all positions above the highest unequal letter cancel exactly in Eq.~\eqref{eq:spectral}.
\end{lemma}
\begin{proof}
If $h\in\ker\chi_i$, the translation $T_h\otimes I_{n_i}$ commutes with
$\diag_e A_{S_{i,\chi_i(e)}}$ and with every translation of the word-block labels in
Eq.~\eqref{eq:Bi}.  It therefore commutes with $\mathcal B_i^t$.
Its diagonal entries at $(e,x_i)$ and $(e+h,x_i)$ agree.  The fibers of
$\chi_i$ are the cosets of its kernel; if $\chi_i=0$, the assertion holds
for every translation.  Equal higher letters include both equal colors and
equal seed-vertex labels, proving the cancellation claim.
\end{proof}

\begin{lemma}[Disjoint binary subspaces]\label{lem:subspaces}
Fix $h_0\ge1$ and let $M=2^{h_0}-1$.  For $k\ge1$, put $d=h_0\lceil\log_2(k+1)/h_0\rceil$.  Whenever $p\ge2d+1$, there are $k$ different $h_0$-dimensional binary subspaces of $\FF_2^p$, all contained in the bit-$0$ zero hyperplane, and any two meeting only in zero.
\end{lemma}
\begin{proof}
The specified $d$ is divisible by $h_0$ and satisfies $2^d>k$.  Identify $\FF_2^{2d}$ with $\FF_{2^d}^2$.  If $K_0=\FF_{2^{h_0}}\subseteq \FF_{2^d}$ and $\lambda_1,\ldots,\lambda_k$ are distinct field elements, use
\[
 V^{(\nu)}=\{(x,\lambda_\nu x):x\in K_0\},\qquad 1\le\nu\le k.
\]
Each subspace has dimension $h_0$, and two distinct subspaces meet only in zero.  Embed the $2d$ coordinates in bit positions $1,\ldots,p-1$.
\end{proof}

For each amplified upper position $j\ge p$, designate a different subspace $\mathcal V_j$ from Lemma~\ref{lem:subspaces} and use all its $M$ nonzero vectors as connection generators at level $i=p+j$.  Extra generators (at most $124$ per position) encode an integer digit.  All proper generator vectors use only the first $p$ label bits, so they fix the seed color at every upper position $j\ge p$.  Thus after division by $q_i$, the matrix~\eqref{eq:Bi} at an amplified upper position has the form
\begin{equation}\label{eq:Zi}
 \mathcal B_i/q_i=\diag_e A_{S_{i,\chi_i(e)}}+\mathcal Z_i\otimes I_{n_i},
\end{equation}
where
\[
 \mathcal Z_i=\frac{D_i}{q_i}I_\Gamma+
 \sum_{s:a_s\le i}\frac{q_{a_s}}{q_i}T_s
\]
is a nonnegative Cayley convolution on $\Gamma$.  At $i=p+j$ with $j\ge p$, every proper displacement lies in $\ker\chi_i$, so the two summands commute.  For a fixed pair of roots $(e,x),(f,y)$, define
\[
 \delta_i=\tau_{S_{i,\chi_i(e)}}(x_i)-\tau_{S_{i,\chi_i(f)}}(y_i),\qquad
 \Delta_i^{(s)}=(A_{S_{i,\chi_i(e)}}^s)_{x_ix_i}-(A_{S_{i,\chi_i(f)}}^s)_{y_iy_i}.
\]
  Write
\[
 \rho=\frac n{n-1}=\frac{125}{124},\qquad
 \kappa=\rho(M+124)+2.
\]
The row-sum estimate is uniform in the word depth. To see this,
let $i_0=2p$ be the first amplified upper position and put
$h_a=\#\{s:a_s=a\}$. Every level $a<i_0$ is an ordinary radix
position and uses at most $n_a-1$ generators. The defining relation
$q_{a+1}=n_aq_a$ therefore telescopes:
\begin{equation}\label{eq:lowerlevelbudget}
 \sum_{a<i_0}h_aq_a
 \le\sum_{a<i_0}(n_a-1)q_a
 =\sum_{a<i_0}(q_{a+1}-q_a)=q_{i_0}-1.
\end{equation}
For $i\ge i_0$, every position from $i_0$ to $i$ has radix
$n=125$ and at most $M+124$ generators. Hence
\begin{equation}\label{eq:highlevelbudget}
 \begin{split}
 \sum_{a\le i}h_a\frac{q_a}{q_i}
 &=\frac1{q_i}\sum_{a<i_0}h_aq_a
   +\sum_{a=i_0}^{i}h_a\frac{q_a}{q_i}\\
 &< 1+(M+124)\sum_{t=0}^{i-i_0}125^{-t}
 <1+\rho(M+124),
 \end{split}
\end{equation}
Since $0\le D_i/q_i<1$,
all row sums of $\mathcal Z_i$ are less than
$\kappa=2+\rho(M+124)$. The same bound applies to every
$i\ge i_0$, independent of $p,k,\ell,r$ and the selected
digit generators.  Each off-diagonal
Cayley coefficient is at most one because every generator displacement
is used once and $q_a\le q_i$; the diagonal coefficient is
$D_i/q_i<1$.  Thus $0\le \mathcal Z_i(e,f)\le1$.  If
$\mu_{i,t}=(\mathcal Z_i^t)_{ee}$, translation invariance makes this
independent of $e$. The generators in $\mathcal V_{i-p}\setminus\{0\}$ have weight one and form a complete graph $K_{M+1}$ on each coset of $\mathcal V_{i-p}$. Since its eigenvalues are $M$ and $-1$, with multiplicities $1$ and $M$, its diagonal return-walk count is $w_t(M)$ below. Entrywise nonnegativity supplies the lower bound.
The entry bound and the row sum bound give
$\mu_{i,t}\le\kappa^{t-1}$ for every $t\ge2$.  Accordingly, for
$t\ge2$,
\begin{equation}\label{eq:mu}
 w_t(M):=\frac{M^t+M(-1)^t}{M+1}\le\mu_{i,t}\le\kappa^{t-1},
\end{equation}
while $\mu_{i,0}=1$ and $0\le\mu_{i,1}\le1$. Since each seed is $d=62$ regular, only powers $s\ge3$ of its adjacency matrix can distinguish seed letters in a diagonal entry.  The binomial theorem gives
\begin{equation}\label{eq:seedmoment}
 \frac{\Delta(\mathcal B_i^g)_{\rm diag}}{q_i^g}
 =\sum_{s=3}^g\binom gs\mu_{i,g-s}
 \Delta_i^{(s)}.
\end{equation}
In particular the first term is $2\binom g3\mu_{i,g-3}\delta_i$.

\begin{lemma}[Uniform high-digit separation]\label{lem:margin}
Fix an integer $g\ge5$. For the mixed-radix construction and generator allocation specified above, there exist constants $M=2^{h_0}-1$, $p_0$, $\lambda_g^->0$ and $\lambda_g^+<\infty$, independent of $r,p,k,\ell$, such that, whenever $p\ge p_0$ and $\istar=p+j$ with $j\ge p$, setting $u=q_{\istar}$ and $\sigma=\operatorname{sgn}(\delta_{\istar})$ gives
\begin{equation}\label{eq:margin}
 \lambda_g^-|\delta_{\istar}|
 \le \sigma\frac{\Delta z_g^R}{u^{g-1}}
 \le\lambda_g^+|\delta_{\istar}|.
\end{equation}
The constants can be chosen so that
\begin{equation}\label{eq:phasewidth}
 4\frac{1+\eta}{1-\eta}\frac{\lambda_g^+}{\lambda_g^-}
 <125^{g-3}.
\end{equation}
\end{lemma}
\begin{proof}
The seeds are $62$-regular and their combined rooted-triangle span is $932$. Hence $|\Delta(A_S^s)_{xx}|\le62^s$ for $s\ge4$ gives
\[
 E_g=\sum_{s=4}^g\binom gs 62^s\kappa^{\max\{g-s-1,0\}},\qquad
 W_g=1864\binom g3\kappa^{g-4}+E_g.
\]
By Lemma~\ref{lem:localcolor}, all positions above $\istar$ cancel exactly.
For every $i\le\istar$, define the normalized (possibly mixed-color)
contribution
\[
 \mathcal R_i=\frac{1}{q_i^g}\bigl[(\mathcal B_i^g)_{(e,x_i),(e,x_i)}
                      -(\mathcal B_i^g)_{(f,y_i),(f,y_i)}\bigr].
\]
The exact decomposition~\eqref{eq:spectral}, without a commutativity
assumption at the lower positions, now reads
\[
 \frac{\Delta z_g^R}{u^{g-1}}
  =\mathcal R_{\istar}
   +\sum_{i<\istar}\left(\frac{q_i}{u}\right)^{g-1}\mathcal R_i.
\]
At the highest amplified position the commuting expansion
\eqref{eq:seedmoment} applies, giving
\[
 \mathcal R_{\istar}=2\binom g3\mu_{\istar,g-3}\delta_{\istar}
                         +\mathcal E_{\istar},\qquad
 |\mathcal E_{\istar}|\le E_g.
\]
The uniform return-walk bounds~\eqref{eq:lowerlevelbudget}--\eqref{eq:mu}
apply at every amplified upper position, independently of the word
depth; hence $|\mathcal R_i|\le W_g$ for earlier amplified positions.
Because their suffix sizes decrease geometrically by $n=125$, their
total normalized contribution is at most
$W_g/(n^{g-1}-1)$.
At an unamplified upper-prefix position, the ordinary-digit generator
budget gives row sum at most $2n$ for $\mathcal B_i/q_i$; at a pad
position it gives row sum at most $2b$. These matrices may mix the two
seed colors, so we use the direct, noncommutative bound
$|\mathcal R_i|\le2(2n_i)^g$ rather than~\eqref{eq:seedmoment}.
The two geometric tails are bounded by
\[
 P_g(p)=\frac{2(2n)^g}{n^{g-1}-1}
 +\frac{2(2b)^g}{(b^{g-1}-1)n^{p(g-1)}}.
\]
Combining the three estimates gives
\[
 \left|\frac{\Delta z_g^R}{u^{g-1}}
   -2\binom g3\mu_{\istar,g-3}\delta_{\istar}\right|
 \le E_g+\frac{W_g}{n^{g-1}-1}+P_g(p).
\]
The highest unequal digit has $|\delta_{\istar}|\ge1$.  Thus the
additive errors can be absorbed into a multiple of $|\delta_{\istar}|$.
The girth $g$ is fixed first, followed by $M$ and $p_0$, before any
choice of $r$ or pair of roots. Take
\begin{equation}\label{eq:lambdas}
\begin{split}
 \lambda_g^-={}&2\binom g3w_{g-3}(M)-E_g
 -\frac{W_g}{n^{g-1}-1}-P_g(p_0),\\
 \lambda_g^+={}&2\binom g3\kappa^{g-4}+E_g
 +\frac{W_g}{n^{g-1}-1}+P_g(p_0).
\end{split}
\end{equation}
We give an explicit sufficient amplification size. Write $B_g=\binom g3$ and take $p_0=3$.
Choose $M=2^{h_0}-1$ so large that
\begin{equation}\label{eq:explicit-amplification}
 M\ge\max\left\{255,\ \frac{16\cdot126^g}{B_g},\
                       \frac{16P_g(3)}{B_g}\right\}.
\end{equation}
This choice depends only on the fixed girth $g$. Since $\rho=125/124$, we have
$\kappa=\rho M+127\le2M$. For every integer $t\ge2$,
\[
 w_t(M)=\frac{M^t+M(-1)^t}{M+1}
 \ge\frac{M-1}{M+1}M^{t-1}\ge\frac{15}{16}M^{t-1}.
\]
For $4\le s\le g$, the exponent $\max\{g-s-1,0\}$ is at most
$g-5$. Consequently,
\[
 E_g\le2^gM^{g-5}\sum_{s=4}^g\binom gs62^s
 \le126^gM^{g-5}\le\frac{B_g}{16}M^{g-4}.
\]
By the third bound in~\eqref{eq:explicit-amplification},
$P_g(p_0)\le(B_g/16)M^{g-4}$. The envelope definition also gives
\[
 W_g\le B_gM^{g-4}\left(1864\cdot2^{g-4}+\frac1{16}\right).
\]
The elementary inequality
\[
 1864\cdot2^{g-4}+\frac1{16}
 <\frac{125^{g-1}-1}{16}\qquad(g\ge5)
\]
holds at $g=5$ and persists because its left side grows by a factor
at most two while the right side grows by more than a factor of two.
Thus $W_g/(125^{g-1}-1)\le(B_g/16)M^{g-4}$.
Substituting these four estimates into~\eqref{eq:lambdas} yields
\begin{equation}\label{eq:explicit-margin}
 \frac{27}{16}B_gM^{g-4}
 \le\lambda_g^-\le\lambda_g^+
 \le\left(2^{g-3}+\frac3{16}\right)B_gM^{g-4}.
\end{equation}
In particular $\lambda_g^->0$ and
$\lambda_g^+/\lambda_g^-<2^{g-3}+1$.
Since $\eta<3/40$ and
\[
 4\frac{43}{37}\bigl(2^{g-3}+1\bigr)<125^{g-3}
 \qquad(g\ge5),
\]
Eq.~\eqref{eq:phasewidth} follows. The last inequality holds at
$g=5$ and propagates to every higher girth: the expression on the
left increases by less than a factor two when $g$ increases by one,
whereas the right side increases by a factor $125$.
The upper and lower bounds in~\eqref{eq:margin} now follow from
\eqref{eq:seedmoment} and the estimates preceding~\eqref{eq:lambdas},
using $|\delta_{\istar}|\ge1$. The parameters $M,p_0$ and
$\lambda_g^\pm$ have all been fixed before $r$ or either root.
The threshold $p_0=3$ controls the spectral error; the disjoint-subspace
allocation has the separate capacity requirement of Lemma~\ref{lem:subspaces},
verified in Lemma~\ref{lem:arith}.
\end{proof}

If $\istar<p$, the highest unequal position is in the pad. If
$\istar=p+j$ with $0\le j<p$, it is in the unamplified upper prefix.
In either case all positions at or below $\istar$ use only ordinary
radix-digit generators. Lemma~\ref{lem:localcolor} cancels every higher
position, and the decomposition~\eqref{eq:spectral} gives
\begin{equation}\label{eq:prefixmoment}
 |\Delta z_g^R|\le C_g^{\mathrm{pre}}u^{g-1},
\end{equation}
where one explicit choice, with $n=125$ and $b=31250$, is
\[
 C_g^{\mathrm{pre}}=
 \frac{2(2n)^g}{1-n^{-(g-1)}}+
 \frac{2(2b)^g}{1-b^{-(g-1)}}.
\]
Indeed, the normalized row sum of $\mathcal B_i/q_i$ at each remaining
position is at most $2n_i$, so $|\mathcal R_i|\le 2(2n_i)^g$.
Summing the upper-prefix and pad contributions as geometric series in
$(q_i/u)^{g-1}$ gives the displayed constant. The largest possible
suffix scale in these cases is $u=b^pn^{p-1}$.

\section{Prescribed degree and phase selection}\label{sec:arithmetic}
For a prescribed $r$, the auxiliary graph must have degree $\omega=N-1-r$ exactly. We first choose a suitable order $N$ and then, only when necessary, choose a suitable representation of that same order. Keeping these two decisions separate is what allows the degree and spectral constraints to be met together.

For $n=125$ and $b=31250=2n^2$, our mixed-radix parameters satisfy
\[
 N=2^\ell b^pn^k=2^{\ell+p}n^{k+2p}.
\]
The identity $b=2n^2$ gives the \emph{scale-preserving change}
\begin{equation}\label{eq:gauge}
 (p,k,\ell)\longmapsto(p+1,k-2,\ell-1).
\end{equation}
It doubles $Q$, halves the number of word blocks and leaves $N$, hence the required deletion degree $\omega=N-1-r$, exactly unchanged. In words, one pad position replaces two upper $n$-positions and one binary label bit. The first irrational rotation selects $N>r$, close enough to $r$ to give an $o(r)$ surplus but with enough room to encode the exact degree. The second moves through these different factorizations of the same $N$ to select a favorable phase, needed only for odd girth. The digit-allocation estimates below hold throughout this second window, so adjusting the phase cannot undo the exact-degree construction. Both choices are made before any roots are compared.

\begin{lemma}[A uniform mesh for irrational rotations]\label{lem:rotation}
Let $\alpha\notin\mathbb Q$ and $\theta>0$. There is an integer $T_0(\alpha,\theta)$ such that, for every $T\ge T_0$, every integer $a$, and every phase $\beta$, the set
\[
 \{\,\beta+(a+t)\alpha\pmod1:0\le t<T\,\}
\]
has maximum circular gap $d_T<\theta$. Thus $d_T\to0$ uniformly over $a$ and $\beta$; also $d_T\ge1/T$.
\end{lemma}
\begin{proof}
Density of the irrational orbit gives a finite initial segment meeting every open arc of length $\theta/2$. Adding points cannot increase its maximal circular gap. Shifting the indices or adding $\beta$ simply rotates the point set, leaving its gaps unchanged. Their sum is one, so the largest of the $T$ gaps is at least $1/T$.
\end{proof}

Two size constraints drive the parameter choice: the word-block size must satisfy $Q=o(N-r)$, and the unamplified prefix must be $o(N^{(g-4)/(g-3)})$. Since $(g-4)/(g-3)\ge1/2$ for every $g\ge5$, it suffices to control this prefix by $O(r^\lambda)$ for some $\lambda<1/2$. Both bounds follow from a common open region of exponent pairs.

\Needspace{9\baselineskip}
\begin{lemma}[Exact-degree and phase selection]\label{lem:arith}
Fix $n=125$, $b=31250$, integers $g\ge5$ and $h_0\ge1$, $M=2^{h_0}-1$, a constant $\gamma>0$, and a nonempty open arc $I$ on the circle $\RR/((g-3)\log n)\mathbb Z$.  For every sufficiently large integer $r$ there are integers $p,k,\ell\ge1$ and a family of distinct generators in~\eqref{eq:R} such that
\begin{enumerate}[label=(\roman*)]
\item $N=2^\ell b^pn^k$ and $r=N-1-\omega$ hold exactly;
\item $p<\ell<k$, $\omega/N\to0$, and for $Q=b^pn^k$ the residual wiring degree $t=\omega-D$ satisfies $t/Q\to\infty$, uniformly in the admissible $p$ window;
\item every amplified upper level uses the $M$ nonzero vectors from its designated complete subspace, with disjoint sets across levels, and at most $124$ additional digit generators. All proper-level generators are supported in label bits $1,\ldots,p-1$, whereas the top-level generators belong to $\{s\in\Gamma:s_0=1\}$;
\item the positive quantities
\begin{equation}\label{eq:Zdef}
 \zeta_j=\gamma\frac{N^{g-4}}{(b^pn^j)^{g-3}},\qquad p\le j<k,
\end{equation}
 have a common logarithmic phase $[\log\zeta_j]$ on this circle, lying in a fixed closed subarc of $I$.
\end{enumerate}
The asymptotic assertions in (ii) are as $r\to\infty$, uniformly over the admissible $p$-window, with $g,h_0,\gamma,I$ fixed.
\end{lemma}
\begin{proof}
Throughout the proof, $M,g,\gamma$ are fixed and $r\to\infty$. Write $k_{\mathrm{tot}}=k+2p$ and $h_{\mathrm{tot}}=\ell+p$; these are invariant under~\eqref{eq:gauge}. We first choose these totals and then vary $p$ in a window on which the degree-allocation estimates remain uniform.

Choose two nondegenerate closed intervals $[x_-,x_+]$ and $[y_-,y_+]$, with $y_->0$, such that every $(x,y)$ in their product satisfies
\begin{equation}\label{eq:feasiblewindows}
 \begin{aligned}
  x\log(2n)-y\log2&>1,\\
  x\log n+2y\log2&<1,\\
  y\log(bn)&<\frac12.
 \end{aligned}
\end{equation}
Such intervals exist: choose $x_0$ strictly between $1/\log(2n)$ and $1/\log n$, then take $y_0>0$ sufficiently small and choose a small closed rectangle around $(x_0,y_0)$. In particular, there are fixed constants $\theta<1$ and $\lambda<1/2$ such that throughout this rectangle
\begin{equation}\label{eq:uniform-exponents}
 x\log n+y\log2<\theta,\qquad y\log(bn)<\lambda.
\end{equation}
Both integer windows $[x_-\log r,x_+\log r]$ and $[y_-\log r,y_+\log r]$ have length $\Theta(\log r)$. The strict inequalities give the required room for both irrational rotations and the degree capacity.

First choose $k_{\mathrm{tot}}$ among the integers in $[x_-\log r,x_+\log r]$, and set
\[
 h_{\mathrm{tot}}=\left\lceil\log_2r-k_{\mathrm{tot}}\log_2n\right\rceil,
 \qquad\xi=h_{\mathrm{tot}}-(\log_2r-k_{\mathrm{tot}}\log_2n).
\]
By Lemma~\ref{lem:rotation} with $\alpha=\log_2n$, the maximal circular gap $d_T$ among the $T=\Theta(\log r)$ successive rotation points tends to zero, uniformly in the position of the integer window, and $d_T\ge1/T$.  The interval $[\sqrt{d_T},3\sqrt{d_T}]$ lies inside $(0,1)$ for large $T$ and has length greater than $d_T$, so it contains a rotation point.  Choose $k_{\mathrm{tot}}$ accordingly.  It follows that
\begin{equation}\label{eq:xi}
 0<\xi\to0,\qquad \xi\ge c(\log r)^{-1/2}.
\end{equation}
Then, with $N=2^{h_{\mathrm{tot}}}n^{k_{\mathrm{tot}}}$,
\begin{equation}\label{eq:near}
 N=r2^\xi>r,\qquad (N-r)/N\to0.
\end{equation}

Next let $p$ range over the integers in $[y_-\log r,y_+\log r]$, setting
\begin{equation}\label{eq:Hk}
 k=k_{\mathrm{tot}}-2p,\qquad \ell=h_{\mathrm{tot}}-p.
\end{equation}
The exact identities $Q=2^pn^{k_{\mathrm{tot}}}$ and $b^pn^p=(bn)^p$, together with~\eqref{eq:uniform-exponents}, yield uniformly over both windows
\begin{equation}\label{eq:uniform-sizes}
 Q\le r^\theta,\qquad b^pn^p\le r^\lambda,
 \qquad\theta<1,\quad\lambda<\frac12.
\end{equation}
Moreover,
\[
 \frac{h_{\mathrm{tot}}}{\log r}
 =\frac{1-(k_{\mathrm{tot}}/\log r)\log n}{\log2}
   +O((\log r)^{-1}).
\]
The second strict inequality in~\eqref{eq:feasiblewindows} gives $2p<h_{\mathrm{tot}}$, and the first gives $h_{\mathrm{tot}}+p<k_{\mathrm{tot}}$, for all sufficiently large $r$, uniformly over the two closed windows. Hence $0<p<\ell<k$. The fixed lower bounds $x_-,y_->0$ also give $p,k=\Theta(\log r)$.

By~\eqref{eq:xi}--\eqref{eq:near}, $N-r\gg r/\sqrt{\log r}$, while~\eqref{eq:uniform-sizes} gives
\[
 \frac{N-r}{Q}\gg\frac{r^{1-\theta}}{\sqrt{\log r}}\longrightarrow\infty.
\]
Since $0\le D<Q$, it follows uniformly that
\begin{equation}\label{eq:tasym}
 t:=N-1-r-D>0,\qquad t/Q\to\infty,\qquad t/N\to0.
\end{equation}

To encode $t$ exactly, proceed along the word positions from low to high.  Start with $t_0=t$. At a pad or low-prefix level put $h_i=t_i\bmod n_i$; at a high level put
\begin{equation}\label{eq:digit}
 h_i=M+((t_i-M)\bmod n_i).
\end{equation}
In both cases set $t_{i+1}=(t_i-h_i)/n_i$. The intermediate $t_i$ need not be nonnegative: Euclidean residues remain well-defined for negative integers, and the telescoping identity holds without that restriction. Nonnegativity is required only for the final top-level coefficient $a_Q$, proved below from $t/Q\to\infty$. At each amplified level select the $M$ prescribed generators plus $h_i-M\in[0,124]$ extra distinct ones.  All other proper levels use $h_i\in[0,n_i-1]$ generators.  Telescoping yields the exact integer decomposition
\begin{equation}\label{eq:telescoping}
 t=\sum_{i<P}h_iq_i+a_QQ,
 \qquad a_Q=t_P.
\end{equation}
The proper-level contribution admits an explicit bound independent of
word length. At an ordinary level $h_i\le n_i-1$, while at an amplified
level $h_i\le(n_i-1)+M$. The mixed-radix telescoping identity gives
\[
 \sum_{i<P}(n_i-1)q_i=Q-1,
 \qquad \sum_{i=2p}^{P-1}q_i=\frac{Q-q_{2p}}{124}.
\]
Consequently, including any negative intermediate Euclidean quotients,
\[
 0\le\sum_{i<P}h_iq_i
 \le Q-1+M\frac{Q-q_{2p}}{124}
 <\left(1+\frac{M}{124}\right)Q.
\]
The exact terminal quotient therefore satisfies
\[
 a_Q=\frac{t-\sum_{i<P}h_iq_i}{Q}
 >\frac tQ-1-\frac{M}{124},\qquad
 a_Q\le\frac tQ.
\]
Since $t/Q\to\infty$ uniformly in the admissible $p$ window, $a_Q$
is a nonnegative integer for every sufficiently large $r$. Also
$a_Q/2^\ell\le t/N\to0$, so the disjoint bit-$0=1$ coset contains
more than $a_Q$ available vectors for top-level $\Pi_P$ generators.

Lemma~\ref{lem:subspaces} gives the disjoint baseline subspaces
provided
\[
 p\ge2h_0\left\lceil\frac{\log_2(k+1)}{h_0}\right\rceil+1,
\]
which holds eventually since $h_0$ is fixed by $g$,
$p=\Theta(\log r)$, and $k=\Theta(\log r)$: the right-hand side is
$O_g(\log\log r)$. This is the allocation condition, separate from
the spectral threshold $p\ge p_0$. There are $k-p$ amplified upper positions,
requiring $M(k-p)$ baseline vectors; all ordinary digits together
require at most $(b-1)p+124k$ further vectors. Thus the number of
used proper displacements is at most
\[
 M(k-p)+(b-1)p+124k
 <2^{p-1}-1
\]
for all sufficiently large $r$: the left side is $O_M(\log r)$,
while the right side grows exponentially in $p$.  Choose the extra
vectors outside all reserved subspaces and all earlier choices, always
with bit $0$ equal to zero. This gives distinct generators with the
required support and realizes~\eqref{eq:degree} at the prescribed $r$.

It remains to select $p$ for the phase.  Fixing $h_{\mathrm{tot}},k_{\mathrm{tot}}$, Eq.~\eqref{eq:Zdef} gives, modulo $(g-3)\log n$,
\begin{equation}\label{eq:rot}
 \log\zeta_j\equiv\log\gamma+(g-4)\log N
 -(g-3)p\log2.
\end{equation}
Every increment of $p$ rotates the phase by $-\log2/\log n$ in units of a period.  This number is irrational, and the $p$ window contains $\Theta(\log r)$ consecutive integers.  Lemma~\ref{lem:rotation} applies uniformly to this window, including its starting phase, which depends on the previously selected $h_{\mathrm{tot}},k_{\mathrm{tot}}$.  Hence some $p$ lies inside a fixed closed subarc of $I$. The
choices~\eqref{eq:Hk} leave $N$ unchanged and preserve all the uniform
estimates in~\eqref{eq:tasym}.

The selected parameters and generators depend only on $r$. Thus the auxiliary graph is fixed before comparing any pair of roots, and it satisfies all four conclusions.
\end{proof}

\paragraph{Phase-free exact degree.}
The construction in the proof of Lemma~\ref{lem:arith}, through the generator allocation in~\eqref{eq:telescoping}, applies to every integer $p$ in the admissible second window. Thus for any fixed $M=2^{h_0}-1$, conclusions (i)--(iii) of that lemma, including $t/Q\to\infty$, hold without choosing a phase arc. We use this phase-free form whenever no cancellation of leading terms is possible.

\begin{lemma}[Exact degree for a fixed alphabet]\label{lem:arith4}
Fix a regular seed alphabet of order $m\ge3$ and degree $d$, where $m$ is not a power of two and $1\le d<m$. For every sufficiently large integer $r$, the layered construction admits parameters $k>\ell\ge1$ and distinct proper- and top-level generators such that
\[
 Q=m^k,\qquad D=d\frac{Q-1}{m-1},\qquad
 N=2^\ell Q,\qquad r=N-1-D-t,
\]
where $t\ge0$, $t/Q\to\infty$ and $t/N\to0$ as $r\to\infty$. At each proper level at most $m-1$ generators are used.
\end{lemma}
\begin{proof}
Choose fixed real numbers $\alpha,\beta$ with
\begin{equation}\label{eq:radixwindow}
 \frac1{\log(2m)}<\alpha<\beta<\frac1{\log m},
 \qquad \theta=\beta\log m<1.
\end{equation}
Let $r\to\infty$ and choose $k$ among the integers in $[\alpha \log r,\beta \log r]$. Then $Q=m^k\le r^\theta$. The number $\log_2m$ is irrational, since a positive integer with a rational base-two logarithm is a power of two. By Lemma~\ref{lem:rotation}, the maximum circular gap $d_T$ of the $T=\Theta(\log r)$ corresponding rotation points tends to zero uniformly in the starting phase, and $d_T\ge1/T$. For sufficiently large $r$, the interval $[\sqrt{d_T},3\sqrt{d_T}]$ is contained in $(0,1)$ and has length greater than $d_T$. Select $k$ so that
\[
 \xi=\left\lceil\log_2(r/Q)\right\rceil-\log_2(r/Q)
 \in[\sqrt{d_T},3\sqrt{d_T}].
\]
Thus $\xi=o(1)$ and $\xi\gg (\log r)^{-1/2}$.

Since $D<Q$, the shift $\log_2(1+(D+1)/r)=O(r^{\theta-1})=o((\log r)^{-1/2})$ leaves this ceiling unchanged. Set
\[
 \ell=\left\lceil\log_2\frac{r+D+1}{Q}\right\rceil,
 \qquad \xi_D=\ell-\log_2\frac{r+D+1}{Q}.
\]
Then $0<\xi_D=o(1)$ and $\xi_D\gg (\log r)^{-1/2}$. Moreover,
\[
 \frac{1-\beta\log m}{\log2}+o(1)
 \le\frac{\ell}{\log r}\le
 \frac{1-\alpha\log m}{\log2}+o(1)<\alpha.
\]
The lower bound is positive, and $k\ge\alpha \log r$, so $0<\ell<k$ eventually. With $N=2^\ell Q$, define
\[
 t=N-1-D-r=(r+D+1)(2^{\xi_D}-1).
\]
It follows that $t/Q\gg r^{1-\theta}/\sqrt{\log r}\to\infty$ and $t/N=1-2^{-\xi_D}\to0$.

Expand $t$ in base $m$:
\[
 t=\sum_{i=0}^{k-1}h_i m^i+a_QQ,
 \qquad 0\le h_i\le m-1,\quad a_Q\ge0.
\]
Assign the $\sum_i h_i=O_m(k)$ proper-level generators to distinct nonzero vectors in $\{s\in\mathbb F_2^\ell:s_0=0\}$ and the $a_Q$ top-level generators to $\{s:s_0=1\}$. The relation $2^{\ell-1}/k\to\infty$ supplies the proper-level capacity. The top-level capacity follows from $a_Q/2^\ell\le t/N\to0$. The degree identity~\eqref{eq:degree} then gives exactly the prescribed $r$.
\end{proof}

For the four-cycle alphabet, we apply this lemma with $m=1250$ and $d=310$.

\Needspace{8\baselineskip}
\section{Proof of the regularity-spectrum classification}\label{sec:proof}
\paragraph{Parameters and quantifiers.}
For each fixed connected non-star $F$ of order at least three, choose the appropriate seed alphabet. In the large-girth case, fix the amplification size $M$, padding depth $p_0$, separating constants $\lambda_g^\pm$, and any necessary safe phase arc; the localization and quotient-error constants are then determined by these fixed choices. Only afterward does the arithmetic construction choose one auxiliary graph for each prescribed degree. The quantifiers are ordered as follows:
\[
 \begin{gathered}
 \forall F\text{ connected and non-star with }|V(F)|\ge3:\\[-2pt]
 \exists r_F\quad\forall r\ge r_F\quad\exists R_r\quad
 \forall v\ne w\in V(R_r).
 \end{gathered}
\]
Each case below compares a uniform leading gap with an error smaller by $\varepsilon+N^{-1}$; at large girth, the signed sum of the two principal terms is controlled first. The threshold and all constants are independent of the root pair.

\begin{proof}[Proof of Theorem~\ref{thm:main}]
Equation~\eqref{eq:star} handles stars. Fix a connected non-star pattern $F$ of order $h$. The proper-level generator budget is at most $b-1$ in the unamplified mixed-radix branches, $\max\{b-1,M+124\}$ in the amplified branch, and $1249$ in the four-cycle branch. Thus the localization and quotient-error constants are uniform in $r$ and the roots. Take $r_F$ large enough for $N\ge2h$, the generator capacities and the estimates below. For each $r\ge r_F$, choose the auxiliary graph once and compare arbitrary distinct roots.

\smallskip\noindent\emph{Case 1: $F$ contains a triangle.}
Apply the exact-degree consequence of Lemma~\ref{lem:arith} with $M=1$ to choose $R$ for the prescribed $r$. For arbitrary distinct roots in this auxiliary graph, Lemma~\ref{lem:leading} and the general quotient bound imply
\[
 a_F\Delta d_F^{\overline R}
 =-6c_3(F)N^{h-3}\Delta\tau_R
 +O_F((\eps+N^{-1})N^{h-3}u^2).
\]
The leading term has magnitude at least $6c_3(F)(1-\eta) N^{h-3}u^2$. Since $\varepsilon\to0$, the error is smaller than this bound for all sufficiently large $r$, uniformly over the roots.

\smallskip\noindent\emph{Case 2: $F$ is a non-star tree.}
Choose $R$ with $M=1$ as in Case~1, and fix arbitrary distinct roots. Such a tree contains a four-vertex path, so $p_4(F)>0$.  There is no uncontracted cyclic support, and the only root-dependent leading quotient type comes from identifying the endpoints of a $P_4$.  Thus
\[
 a_F\Delta d_F^{\overline R}
 =8p_4(F)N^{h-4}\Delta\tau_R
 +O_F((\eps+N^{-1})N^{h-4}u^2).
\]
Its leading term has magnitude at least $8p_4(F)(1-\eta) N^{h-4}u^2$, so the uniform error is smaller for all sufficiently large $r$.

\smallskip\noindent\emph{Case 3: $F$ has no triangle but has a four-cycle.}
Use Lemma~\ref{lem:arith4} and the $U/V$ seed system to choose $R$. For arbitrary distinct roots, Lemmas~\ref{lem:c4} and~\ref{lem:c4main} give
\[
 a_F\Delta d_F^{\overline R}
 =8c_4(F)N^{h-4}\Delta c_4^R
 +O_F((\eps+N^{-1})N^{h-4}u^3).
\]
The fixed positive gap in~\eqref{eq:c4inv} makes the leading magnitude a positive multiple of $N^{h-4}u^3$. The error is smaller for all sufficiently large $r$, uniformly over the roots.

\smallskip\noindent\emph{Case 4: $F$ has finite girth $g\ge5$.}
Fix $M,p_0,\lambda_g^-,\lambda_g^+$ as in Lemma~\ref{lem:margin}. For even $g$, choose $R$ for the prescribed $r$ by the phase-free exact-degree construction following Lemma~\ref{lem:arith}: the two leading contributions will reinforce one another. For odd $g$, put
\[
 \gamma_F=\frac{8p_4(F)}{gc_g(F)}>0
\]
and choose a fixed nonempty open arc $I$ disjoint from the image of
\[
 \left[\log\frac{\lambda_g^-}{2(1+\eta)},\,
       \log\frac{2\lambda_g^+}{1-\eta}\right]
\]
in $\RR/((g-3)\log n)\mathbb Z$. Such an arc exists by Eq.~\eqref{eq:phasewidth}. Apply Lemma~\ref{lem:arith} with $\gamma=\gamma_F$ and this arc. In either parity, the construction fixes one $R$ before comparing roots. Increase $r_F$ so that $p\ge p_0$ and the uniform error bounds below hold.

Fix arbitrary distinct roots $v,w$ of this auxiliary graph. Write $u=q_{\istar}$ and set
\[
\begin{aligned}
 S_P&=N^{h-4}u^2,& S_C&=N^{h-g}u^{g-1},\\
 T_P&=8p_4(F)N^{h-4}\Delta\tau_R,&
 T_C&=(-1)^g gc_g(F)N^{h-g}\Delta z_g^R.
\end{aligned}
\]
These are the scales and signed contributions in~\eqref{eq:principal}. The ratio $S_C/S_P$ governs the comparison. At low word positions the cycle scale becomes negligible, leaving the contracted-path term. At amplified positions both scales may matter: for even $g$ their signed contributions reinforce, while for odd $g$ the chosen phase keeps $\zeta_j=\gamma_F S_P/S_C$ away from the cancellation range. The estimates below make all three cases uniform over the roots.

If $\istar$ lies in the pad or in the unamplified upper prefix, then $u\le b^pn^p\le r^\lambda$ by~\eqref{eq:uniform-sizes}, where $\lambda<1/2$. As $N>r$ and $g\ge5$ is fixed,
\[
 \vartheta:=\frac{S_C}{S_P}=\frac{u^{g-3}}{N^{g-4}}
 \le r^{\lambda(g-3)-(g-4)}\longrightarrow0,
\]
since $\lambda(g-3)<g-4$ for every $g\ge5$.
Equations~\eqref{eq:trigap}, \eqref{eq:prefixmoment} and \eqref{eq:principal} give
\[
 |T_P|\ge8p_4(F)(1-\eta) S_P,\qquad |T_C|\le C_FS_C,
\]
and therefore
\begin{equation}\label{eq:lowabsorb}
 \frac{|T_C|+|E_F|}{|T_P|}
 \le C'_F\bigl[\vartheta+
 (\varepsilon+N^{-1})(1+\vartheta)\bigr]\longrightarrow0.
\end{equation}
Thus the total root difference is nonzero in this case.

Now suppose $\istar=p+j$ with $j\ge p$, so $u=b^pn^j$, and put $\sigma=\operatorname{sgn}(\delta_{\istar})$. Only when $g$ is odd do we also write
\[
 \zeta_j=\gamma_F\frac{N^{g-4}}{u^{g-3}}.
\]
By Lemma~\ref{lem:margin} and the triangle-digit expansion,
\begin{equation}\label{eq:twointerval}
 \begin{split}
 (1-\eta)u^2|\delta_{\istar}|&\le \sigma\Delta\tau_R\le(1+\eta)u^2|\delta_{\istar}|,\\
 \lambda_g^-u^{g-1}|\delta_{\istar}|&\le\sigma\Delta z_g^R\le\lambda_g^+u^{g-1}|\delta_{\istar}|.
 \end{split}
\end{equation}
Define the uniform separation margin
\[
 m_F=\tfrac13\min\{8p_4(F)(1-\eta),\,gc_g(F)\lambda_g^-\}>0.
\]
For even $g$, the two contributions have the same sign $\sigma$. Hence
\begin{equation}\label{eq:evenabsorb}
\begin{split}
 |T_P+T_C|&=\sigma(T_P+T_C)\\
 &\ge |\delta_{\istar}|\bigl[8p_4(F)(1-\eta)S_P
                  +gc_g(F)\lambda_g^-S_C\bigr]\\
 &\ge m_F(S_P+S_C).
\end{split}
\end{equation}

For odd $g$, the phase arc chosen before the roots ensures that every amplified position satisfies
\begin{equation}\label{eq:safes}
 \zeta_j\le\frac{\lambda_g^-}{2(1+\eta)}
 \quad\text{or}\quad
 \zeta_j\ge\frac{2\lambda_g^+}{1-\eta}.
\end{equation}
In the first range $|T_C|\ge2|T_P|$, and in the second $|T_P|\ge2|T_C|$. Since the signs are opposite,
\[
 |T_P+T_C|=\bigl||T_P|-|T_C|\bigr|
 \ge\tfrac13(|T_P|+|T_C|)
 \ge m_F(S_P+S_C).
\]
For sufficiently large $r$, the error in~\eqref{eq:principal} is less than half this margin because $C_F(\varepsilon+N^{-1})<m_F/2$. Consequently the same graph satisfies, simultaneously for every high-position pair,
\[
 |a_F\Delta d_F^{\overline R}|
 \ge |T_P+T_C|-|E_F|
 \ge \tfrac12 m_F(S_P+S_C)>0.
\]
The threshold needed for this estimate depends on $F$ and its fixed amplification constants, never on the pair of roots.

These cases exhaust the connected non-star graphs.  Finally, for every branch $\omega/N\to0$ and therefore
\[
 \frac{r}{N}=1-\frac{1+\omega}{N}\longrightarrow1.
\]
For large $r$ we have $r>(N-2)/2$; hence every two vertices of
$G=\overline R$ have distance at most two, by the elementary
common-neighbor argument. The graph $R$ has an edge, so $G$ is not
complete and its diameter is exactly two.

In every case, the uniform separation estimates hold for the single graph chosen at each $r$. Together with the exact degree, asymptotic order and diameter estimates, this proves the theorem.
\end{proof}

\section{Quantitative extremal bounds for triangles}\label{sec:quantitative}
For triangles, word encoding and shell invariance suffice without the quotient localization and spectral amplification needed for general patterns. Two equal-degree seeds of order $29$ give word blocks of size $29^k$ and $2^k$ color words, hence the total scale $58^k$ and the surplus exponent $\log29/\log58$. The finite certificate checks the seed properties; the word and degree constructions apply at arbitrary depth.

\subsection{A complementary 29-vertex alphabet}
Let $S_0$ be the graph on $\{0,\ldots,28\}$ specified by the adjacency masks in Appendix~\ref{app:qt-certificate}, and put $S_1=\overline{S_0}$.
\begin{lemma}[The quantitative seed pair]\label{lem:qt-seed}
Both $S_0$ and $S_1$ are $14$-regular, and
\begin{equation}\label{eq:qt-seeds}
 \begin{split}
 \{\tau_{S_0}(v):v\in V(S_0)\}&=\{43,44,\ldots,71\},\\
 \{\tau_{S_1}(v):v\in V(S_1)\}&=\{13,14,\ldots,41\}.
 \end{split}
\end{equation}
The $58$ counts, taken jointly, are all distinct and have total span $58$.
\end{lemma}
\begin{proof}
The verifier \texttt{verify\_triangle29.py}, supplied with the finite certificates, reconstructs $S_0$ from the masks in Appendix~\ref{app:qt-certificate}, checks symmetry, zero diagonal and degree $14$, and enumerates every unordered vertex triple. It gives the first set in~\eqref{eq:qt-seeds}. Equation~\eqref{eq:triangle-complement}, with $N=29$ and $d=14$, gives $\tau_{S_1}(v)=84-\tau_{S_0}(v)$, yielding the second set and joint injectivity.
\end{proof}

For a small illustration, $\overline{S_0\mathbin{\dot\cup}S_1}$ is a $43$-regular graph on $58$ vertices, with rooted triangle counts
\[
 \{622,\ldots,650\}\ \cup\ \{652,\ldots,680\}.
\]
Adding the matching between equally labelled vertices of $S_0$ and $S_1$ before taking the complement gives a $42$-regular triangle-distinct graph of the same order. That matching creates no triangles in the disjoint union. The large-degree construction repeats the same principle at several suffix scales.

\subsection{Word encoding}
Throughout this section put $q=29$. For a word $c=(c_0,\ldots,c_{k-1})\in\{0,1\}^k$, let
\[
 W_c=S_{c_{k-1}}[\cdots[S_{c_0}]\cdots].
\]
Write $q_i=q^i$, $Q=q^k$ and
\begin{equation}\label{eq:qt-worddegree}
 D=14\sum_{i=0}^{k-1}q^i=\frac{Q-1}{2}.
\end{equation}
Every suffix of order $q_i$ has degree $D_i=(q_i-1)/2$. Iterating~\eqref{eq:trianglelex} gives
\begin{equation}\label{eq:qt-digits}
 \tau_{W_c}(x)=C_k+\sum_{i=0}^{k-1}q_i^2\tau_{S_{c_i}}(x_i),
\end{equation}
where $C_k$ is independent of the word and root.
\begin{lemma}[Joint separation of words and roots]\label{lem:qt-words}
The $2^kQ$ rooted triangle counts of all the $W_c$, taken jointly across their color words, are pairwise distinct.
\end{lemma}
\begin{proof}
At the highest differing word-root coordinate $i$, the seed counts differ by a nonzero integer. Its weighted contribution has absolute value at least $q_i^2$. The lower-coordinate contribution has absolute value at most
\[
 58\sum_{j<i}q_j^2=\frac{58(q_i^2-1)}{29^2-1}<q_i^2,
\]
so it cannot cancel the leading difference.
\end{proof}

For $L\le2^k$, assign $L$ distinct color words to $\mathbb Z/L\mathbb Z$ and use the shell matrices from~\eqref{eq:shells}, now with $q_i=29^i$ and $P=k$. Lemma~\ref{lem:shell-invariance} shows that every symmetric shell connection preserves the joint triangle-count separation. Its row degree is independent of the root. It remains to choose these row degrees so that the complement has the prescribed regularity.

\subsection{Simultaneous interpolation of all correction digits}
The next observation applies to any fixed radix. A single inverse-displacement pair can carry shells from several levels, because their supports are disjoint.
\begin{lemma}[Constant inverse-pair budget]\label{lem:qt-shell-digits}
Let $q\ge2$, $k\ge1$, and let $E_0,\ldots,E_k$ be the suffix shells for $Q=q^k$. For every integer $0\le H<Q$ there are at most $q-1$ shell connections $K_1,\ldots,K_{q-1}$ such that
\begin{equation}\label{eq:qt-shell-sum}
 \sum_{s=1}^{q-1}\operatorname{rowSum}(K_s)=H.
\end{equation}
\end{lemma}
\begin{proof}
Write the shell-digit expansion
\[
 H=h_0+(q-1)\sum_{j=1}^{k}h_jq^{j-1},
 \qquad 0\le h_0\le q-2,\quad 0\le h_j\le q-1\ (1\le j\le k).
\]
Indeed, $h_0$ is the residue of $H$ modulo $q-1$, and the quotient
$(H-h_0)/(q-1)$ is less than $q^k$, so it has at most $k$ base-$q$
digits. For $1\le s\le q-1$, put
\begin{equation}\label{eq:qt-shell-assignment}
 K_s=\sum_{j=0}^{k}\mathbf1_{\{s\le h_j\}}E_j.
\end{equation}
The shell degrees are $1,(q-1),(q-1)q,\ldots,(q-1)q^{k-1}$. Summing the row degrees over $s$ yields
\[
 \sum_{j=0}^{k}h_j\operatorname{rowSum}(E_j)
 =h_0+(q-1)\sum_{j=1}^{k}h_jq^{j-1}=H.
\]
Each $K_s$ is zero--one, since its shells have disjoint supports. Different $K_s$ will be attached to different displacements, so their edges in the full graph are also disjoint.
\end{proof}

\subsection{All sufficiently large regularities}
\begin{proof}[Proof of Theorem~\ref{thm:quantitative}]
Define
\[
 \mathfrak r_k=58^k-\frac{29^k+1}{2}.
\]
This sequence is strictly increasing and unbounded. Given $r\ge r_\triangle=\mathfrak r_{11}+1$, choose the unique $k\ge12$ for which $\mathfrak r_{k-1}<r\le\mathfrak r_k$. Put $Q=29^k$, $D=(Q-1)/2$, and first use the even label count
\begin{equation}\label{eq:qt-parameters}
 L=2\left\lceil\frac{r+D+1}{2Q}\right\rceil,
 \qquad N=LQ,\qquad t=N-1-r-D.
\end{equation}
The upper inequality $r\le\mathfrak r_k=2^kQ-D-1$ ensures $L\le2^k$, while rounding upward to an even multiple gives $0\le t<2Q$. The lower endpoint $r>\mathfrak r_{k-1}$ gives
\begin{equation}\label{eq:qt-capacity}
 \frac{r+D+1}{Q}>
 \frac{2^{k-1}+14}{29}\ge\frac{2062}{29}>71,
\end{equation}
and hence $L\ge72$. The $28$ inverse pairs $\{s,-s\}$, $1\le s\le28$, are therefore mutually disjoint and avoid the involution $L/2$.

Write $t=\nu+2H$, with $\nu\in\{0,1\}$ and $0\le H<Q$. Apply Lemma~\ref{lem:qt-shell-digits} with $q=29$, assigning $K_s$ to both displacements $s$ and $-s$. When $\nu=1$, assign $E_0=I_Q$ to $L/2$. The resulting cross-block degree is exactly
\[
 2\sum_{s=1}^{28}\operatorname{rowSum}(K_s)+\nu=t.
\]
The auxiliary graph has degree $\omega=D+t=N-1-r$. Lemmas~\ref{lem:qt-words} and~\ref{lem:shell-invariance} make it triangle-distinct. By~\eqref{eq:triangle-complement}, $G_r=\overline R$ remains triangle-distinct and is exactly $r$-regular. Equation~\eqref{eq:p4exact} gives $P_4$-irregularity.

For the order estimate,
\[
 0<N-r=1+D+t<\tfrac52Q.
\]
Since $k\ge12$,
\[
 \mathfrak r_{k-1}\ge(1-2^{-11})58^{k-1},\qquad
 Q<29\left(\frac{r}{1-2^{-11}}\right)^{\beta_\triangle}.
\]
Using $0<\beta_\triangle<1$, we obtain
\begin{equation}\label{eq:qt-73}
 N-r<\frac{72.5}{(1-2^{-11})^{\beta_\triangle}}r^{\beta_\triangle}
 <\frac{72.5}{1-2^{-11}}r^{\beta_\triangle}<73r^{\beta_\triangle}.
\end{equation}
Also $\omega<\tfrac52Q$ and $N\ge72Q>2\omega$. Thus $2r>N-2$, and any two nonadjacent vertices of $G_r$ have a common neighbor. Since $\omega>0$, the graph is not complete and its diameter is exactly two.

For even $r$, take instead
\begin{equation}\label{eq:qt-even-labels}
 L=\left\lceil\frac{r+D+1}{Q}\right\rceil.
\end{equation}
Now $0\le t=LQ-1-r-D<Q$. As $Q$ is odd and $D$ is even, $t$ is even when $L$ is odd and odd when $L$ is even. Thus an involution is needed only when it exists. Inequality~\eqref{eq:qt-capacity} still gives $L\ge72$, and $L\le2^k$ still holds. There are at least $28$ inverse pairs in either parity, so the same assignment applies. Finally,
\[
 N-r=1+D+t<\tfrac32Q
 <\frac{43.5}{(1-2^{-11})^{\beta_\triangle}}r^{\beta_\triangle}
 <44r^{\beta_\triangle}.
\]
The separation and diameter arguments are unchanged.
\end{proof}

\section{Finite seed lifting and exponent constraints}\label{sec:alphabet}
The shell interpolation applies to any equal-degree seed alphabet with jointly distinct rooted triangle counts. We give the general lifting statement and then a necessary restriction on such alphabets.

\subsection{Lifting an equal-degree seed collection}
\begin{theorem}[Finite seed lifting]\label{thm:qt-lifting}
Suppose $a\ge2$ regular graphs $S_1,\ldots,S_a$ have the same order $q\ge3$ and degree $d>0$, and their $aq$ rooted triangle counts are jointly pairwise distinct. Then every sufficiently large integer $r$ admits a connected, diameter-two, $r$-regular triangle-distinct graph of order
\begin{equation}\label{eq:qt-lifting}
 N=r+O_{a,q,d}\left(r^{\log q/\log(aq)}\right).
\end{equation}
Its degree correction uses at most $q-1$ inverse-displacement pairs, independently of word length. Odd values of $a$ are allowed.
\end{theorem}
\begin{proof}
Form the $a^k$ possible color words of length $k$, and put $Q=q^k$ and $D=d(Q-1)/(q-1)$. The rooted-count identity~\eqref{eq:qt-digits} holds for every word. All seed counts lie between $0$ and $\binom d2<q^2-1$; hence the highest differing digit dominates the geometric tail of lower digits, as in Lemma~\ref{lem:qt-words}. We shall select an even number of these color words for the cyclic label group.

Set $A_k=2\lfloor a^k/2\rfloor$ and $\mathfrak r_k=A_kq^k-d(q^k-1)/(q-1)-1$. The sequence is strictly increasing and unbounded. Indeed, $A_{k+1}\ge A_k\ge2$ and
\[
 \mathfrak r_{k+1}-\mathfrak r_k
 =q^k(qA_{k+1}-A_k-d)>0,
\]
since $A_k\ge2$ and $d\le q-1$. Choose $k$ with $\mathfrak r_{k-1}<r\le\mathfrak r_k$, and select $L,N,t$ by~\eqref{eq:qt-parameters}. Then $L\le A_k$, $0\le t<2Q$, and $t=\nu+2H$ with $\nu\in\{0,1\}$ and $0\le H<Q$.

Lemma~\ref{lem:qt-shell-digits} encodes $H$ with $q-1$ inverse pairs. The inequality
\[
 \frac{r+D+1}{Q}>\frac{A_{k-1}+d}{q}
\]
shows that the sufficient capacity condition
\begin{equation}\label{eq:qt-general-capacity}
 A_{k-1}+d\ge2q(q-1)
\end{equation}
gives $L\ge2q$ and hence $L/2-1\ge q-1$. The required pair budget, and the right-hand side of~\eqref{eq:qt-general-capacity}, do not grow with $k$. The condition holds for all sufficiently large $k$ because $a\ge2$. If $\nu=1$, the involution contributes $I_Q$.

Lemma~\ref{lem:shell-invariance} applies with the common suffix degrees $D_i=d(q^i-1)/(q-1)$. The shell connections therefore retain the joint separation of the chosen words, while their row degrees give $\omega=D+t=N-1-r$ exactly. Complementation preserves injectivity. Moreover, $D<Q$, $t<2Q$ and $L\ge2q\ge6$ imply $N>2\omega$; since $\omega>0$, the complement is connected and has diameter exactly two. Finally $N-r=O(Q)$ and $\mathfrak r_{k-1}\asymp(aq)^{k-1}$, whence $Q=O(r^{\log q/\log(aq)})$.
\end{proof}

\subsection{A necessary seed-size bound}
We apply the triangle-degree estimate of Hak, Kozerenko and Serdiuk~\cite[Lemma~3.5]{Hak2025} to the seeds in a common alphabet. Its short proof is recalled for completeness.
\begin{lemma}[Known triangle-degree bound]\label{lem:qt-twins}
If a $d$-regular graph $H$ on at least two vertices is triangle-distinct, then
\begin{equation}\label{eq:qt-twin-bound}
 \tau_H(v)\le\binom d2-\left\lceil\frac d2\right\rceil
 \qquad(v\in V(H)).
\end{equation}
\end{lemma}
\begin{proof}
Let $h=\binom d2-\tau_H(v)$ count the nonedges inside $N_H(v)$. If $h<\lceil d/2\rceil$, fewer than $d$ neighbors are incident with those nonedges. Some $u\in N_H(v)$ is adjacent to every other neighbor of $v$. Its $d$ neighbors are precisely $v$ and $N_H(v)\setminus\{u\}$, so $N_H[u]=N_H[v]$. Swapping $u$ and $v$ is an automorphism, contradicting distinct triangle counts.
\end{proof}

\begin{theorem}[Size of an equal-degree triangle alphabet]\label{thm:qt-alphabet}
Under the seed hypotheses of Theorem~\ref{thm:qt-lifting}, set $s=\min\{d,q-1-d\}$. Then
\begin{equation}\label{eq:qt-alphabet-capacity}
 aq\le\binom s2-\left\lceil\frac s2\right\rceil+1,
 \qquad q\ge8a+8.
\end{equation}
\end{theorem}
\begin{proof}
Put $w=q-1-d$ and $B_\triangle(j)=\binom j2-\lceil j/2\rceil$. Equation~\eqref{eq:triangle-complement} gives
\[
 \tau_{\overline{S_i}}(v)=C-\tau_{S_i}(v),\qquad
 C=\binom{q-1}{2}-\tfrac32dw,
\]
with the same constant for every seed. Their complements are thus jointly triangle-distinct. By Lemma~\ref{lem:qt-twins}, all $aq$ counts lie in
\[
 [0,B_\triangle(d)]\cap[C-B_\triangle(w),C],
\]
which contains at most $B_\triangle(s)+1$ integers. The seed hypotheses force $d,w\ge3$, and $B_\triangle$ is increasing in that range. This proves the first inequality.

Since $q\ge2s+1$,
\[
 a(2s+1)\le aq\le B_\triangle(s)+1\le\tfrac12s(s-2)+1,
\]
so $0\le\tfrac12s(s-4a-2)+(1-a)$. If $s\le4a+2$, the right-hand side is negative. Hence $s\ge4a+3$ and $q\ge8a+7$. Equality would force the odd order $q=8a+7$ and odd degree $d=w=4a+3$, contradicting the handshaking lemma. Therefore $q\ge8a+8$.
\end{proof}

In particular, a two-seed alphabet requires $q\ge24$, and a three-seed alphabet requires $q\ge32$. At order $29$ the alphabet of Lemma~\ref{lem:qt-seed} attains the maximum possible number of colors. For fixed $a>1$, $\log q/\log(aq)$ increases with $q$, so this lifting mechanism satisfies
\begin{equation}\label{eq:qt-exponent-barrier}
 \frac{\log q}{\log(aq)}\ge
 \frac{\log(8a+8)}{\log(a(8a+8))}.
\end{equation}
For $a=2$, the right-hand side is $\log24/\log48=0.820948\ldots$, whereas our exponent is $\beta_\triangle=0.829292720\ldots$. The bound concerns the equal-degree seed-lifting construction.

\subsection{Comparison with the necessary order surplus}
The triangle-degree restrictions of Berikkyzy et al.~\cite[Theorem~5]{Berikkyzy2024} and Hak, Kozerenko and Serdiuk~\cite[Lemma~3.5 and Corollary~3.6]{Hak2025} yield a lower bound on the order surplus. If $G$ is $r$-regular and triangle-distinct, with order $N$, its complement has degree $w=N-1-r$ and is also triangle-distinct. Lemma~\ref{lem:qt-twins} implies
\[
 N-1=r+w\le\binom w2-\left\lceil\frac w2\right\rceil
 \le\frac{w(w-2)}2.
\]
Solving for $w$ gives
\begin{equation}\label{eq:qt-lower}
 N-r\ge\left\lceil3+\sqrt{2r+4}\right\rceil.
\end{equation}
The square-root lower bound and Theorem~\ref{thm:quantitative} thus leave an exponent gap between the necessary and constructive order surpluses.

\section{Consequences and open questions}\label{sec:conclusion}
The classification theorem leaves two extremal questions. For a fixed connected non-star $F$, determine the finite exceptional set $\mathbb Z_{\ge0}\setminus\Spec(F)$, or bound its largest element effectively. For triangles, narrow the gap between the square-root order obstruction~\eqref{eq:qt-lower} and the constructive exponent $\log29/\log58$ in Theorem~\ref{thm:quantitative}. The seed-size restriction in Theorem~\ref{thm:qt-alphabet} identifies one limitation of the present lifting method.

\appendix
\section{Finite inputs and verification protocol}\label{app:certificates}
Lemma~\ref{lem:seeds} uses the joint injectivity of the rooted triangle counts for $A,B$ and of the rooted four-cycle counts for $U,V$. Appendix~\ref{app:embedded} gives their constructions, which also establish regularity and edge-disjointness; the injectivity for $H_0$ follows from~\eqref{eq:trianglelex}. The separate $29$-vertex pair for the quantitative construction is specified in Appendix~\ref{app:qt-certificate}.

Appendix~\ref{app:embedded} contains the finite reconstruction rules and the block fourth-moment formula. The deterministic verifier \texttt{verify\_general\_seeds.py} checks degrees, complementarity, disjointness and rooted counts, including agreement with direct common-neighbor enumeration. The separate \texttt{verify\_triangle29.py} checks the 29-row adjacency certificate. Both scripts use only the Python standard library and are supplied with the supplementary material.

\paragraph{Exact-rational examples for the spectral margin.}
For $M=2^{30}-1$, $p_0=3$, $n=125$ and $b=31250$, substitution
into~\eqref{eq:mu} and~\eqref{eq:lambdas}, with
$L_g=2\binom g3 M^{g-4}$, gives the exact bounds
\[
 \frac{19}{20}L_g<\lambda_g^-<\lambda_g^+<\frac{53}{50}L_g,
 \qquad g\in\{5,6,7\}.
\]
Clearing the positive denominators in the definitions of
$E_g,W_g,P_g,\kappa$ and $w_{g-3}(M)$ gives these inequalities. Since $\eta<3/40$,
they imply the stronger explicit phase test
\[
 4\frac{1+\eta}{1-\eta}\frac{\lambda_g^+}{\lambda_g^-}
 <4\frac{43}{37}\frac{53/50}{19/20}<6<125^{g-3}.
\]
The supplementary Python verifier checks these rational cross-products
using exact \texttt{Fraction} arithmetic.

\section{Fourth-order connection algebra}\label{app:fourth}
Here we record the matrix identities underlying Lemma~\ref{lem:c4}. At a fixed connection level $a$, write $W_e=X_e[Y_e]$, with $|V(Y_e)|=q=q_a$ and $d(Y_e)=D_a$. The associated kernel is $\Pi_a=I\otimes J_q$ from~\eqref{eq:proj}. All suffix graphs at this level have the same order and degree. Directly from
\[
 A_{W_e}=A_{X_e}\otimes J_q+I\otimes A_{Y_e},
 \qquad J_qA_{Y_e}=D_aJ_q,
\]
we obtain the diagonal identities
\begin{align}
 (\Pi_a A_{W_e})_{xx}&=D_a,\\
 (\Pi_a A_{W_e}\Pi_a)_{xx}&=qD_a,\\
 (\Pi_a A_{W_e}^2\Pi_a)_{xx}&=q^3d(X_e)+qD_a^2.\label{eq:PAP2}
\end{align}
The only mixed term in a fourth-order closed walk that retains information about the two seed colors is
\begin{equation}\label{eq:PABP}
 (\Pi_a A_{W_e}A_{W_{e+s}}\Pi_a)_{xx}
 =q^3(A_{X_e}A_{X_{e+s}})_{x_+,x_+}+qD_a^2.
\end{equation}
For clarity, write $H=A_{X_e}$, $H'=A_{X_{e+s}}$ and $J=J_q$.  Multiplying
$A_{W_e}=H\otimes J+I\otimes A_{Y_e}$ and
$A_{W_{e+s}}=H'\otimes J+I\otimes A_{Y_{e+s}}$
using $J A_{Y_e}=A_{Y_e}J=D_aJ$ gives the full identity
\[
 \Pi_a A_{W_e}A_{W_{e+s}}\Pi_a
 =\bigl(q^3HH'+q^2D_a(H+H')+qD_a^2I\bigr)\otimes J.
\]
Both $H$ and $H'$ have zero diagonal. Consequently the mixed $q^2D_a$ term vanishes on the diagonal, leaving exactly the $q^3$ and $qD_a^2$ contributions in~\eqref{eq:PABP}. The first is a same-label common-neighbor count in the high word graphs. The lexicographic recursion, for equal-order equal-degree seed pairs, is
\begin{equation}\label{eq:commonneighbors}
 \operatorname{cn}_{X[Y],X'[Y']}((x,y))
 =|V(Y)|\operatorname{cn}_{X,X'}(x)+\operatorname{cn}_{Y,Y'}(y),
\end{equation}
where $\operatorname{cn}_{H,H'}(v)=\sum_zA_H(v,z)A_{H'}(v,z)$.
\begin{lemma}[Mixed-color common-neighbor recursion]\label{lem:mixed-cn}
Suppose at position $i$ the two seed graphs $S_{i,0}$ and $S_{i,1}$
have common order $n_i$ and degree $d_i$, and their edge sets are
disjoint. For the associated lexicographic words $W_e,W_f$ on the
same labeled vertex set, every word vertex $x$ satisfies
\[
 \operatorname{cn}_{W_e,W_f}(x)
 =\sum_{i=0}^{P-1}q_i d_i
       \mathbf1_{\{\chi_i(e)=\chi_i(f)\}}.
\]
In particular, for $f=e+s$ in a binary label group, this expression
depends on $s$ and the sequence of levels, but not on $e$ or $x$.
\end{lemma}
\begin{proof}
For two equal-order lexicographic products $X[Y]$ and $X'[Y']$,
expand their adjacency matrices:
\[
 (A_X\otimes J_q+I\otimes A_Y)
 (A_{X'}\otimes J_q+I\otimes A_{Y'}).
\]
The cross terms contain $A_X$ or $A_{X'}$ on the outer diagonal and
therefore vanish there, since all graphs are loopless. The other two
terms give~\eqref{eq:commonneighbors}. If the colors coincide at
position $i$, the two seed rows have scalar product $d_i$;
if they differ, their scalar product is zero by edge disjointness.
Iterating the recursion from the lowest word coordinate proves the
formula. Finally $\chi_i(e+s)=\chi_i(e)+\chi_i(s)$, so the indicators
are determined by $s$ alone.
\end{proof}
For the $U/V$ alphabet all $d_i=310$ and the stated edge-disjointness
holds. Lemma~\ref{lem:mixed-cn} proves that~\eqref{eq:PABP} depends
only on the displacement $s$ and the connection level, not on the
rooted letter. To see that this exhausts the fourth
power, write $\mathsf B=\bigoplus_e A_{W_e}$ and let $\mathsf L$ denote the connection matrix.  A fourth-order diagonal word with one $\mathsf L$ cannot return to
its starting word block.  If it has exactly two $\mathsf L$ factors, their
nonzero displacements must be equal in $\FF_2^\ell$; since each
displacement is assigned only one level, the factors have the same
connection level $a$.  By suffix regularity,
$A_{W_e}\Pi_a=\Pi_a A_{W_e}$ for every word block $e$.
Thus the words with two consecutive $\mathsf L$ factors reduce to
$\Pi_a A_{W_e}^2\Pi_a$ or its cyclic reordering, whose diagonal
is constant by~\eqref{eq:PAP2}; the alternating words reduce to
\eqref{eq:PABP}.  Three $\mathsf L$ factors may have different levels, but
return to their original word block only when their three displacements
sum to zero.  Their kernel product is a scalar multiple of
$\Pi_{\max(a,b,c)}$ by~\eqref{eq:projrule}; after commuting with
the unique $\mathsf B$ factor, the diagonal becomes a scalar multiple of
$(\Pi_d A_{W_e})_{xx}=D_d$, again independent of the root.
Four $\mathsf L$ factors contain no $\mathsf B$ and their diagonal, when their
XOR displacement is zero, is constant by~\eqref{eq:projrule}.
This proves root independence of every term in
$\diag((\mathsf B+\mathsf L)^4-\mathsf B^4)$, for every choice of distinct generator displacements.

\section{Cyclic quotient losses and source-root multiplicities}\label{app:quotient-ledger}
We record the quotient-loss classification used in Lemma~\ref{lem:leading}
for an arbitrary fixed source pattern $F$.
For any admissible independent-block quotient of an edge support, use
$s,m,c,d=s-m$, the rooted cyclic girth $q$, and
$A_q=m-c-q+1\ge0$ as in Lemma~\ref{lem:quotient}. A rooted tree
component cancels exactly, even if another quotient component has
a cycle. For the remaining terms the ratios to the candidate
leading scales are:
\begin{center}
\small
\begin{tabular}{@{}p{0.28\linewidth}p{0.27\linewidth}p{0.35\linewidth}@{}}
\toprule
Source and quotient class & Ratio bound & Zero-loss or suppression test \\
\midrule
Any $F$, rooted cyclic quotient, triangle scale
 &$N^{-d}\varepsilon^{q-3+A_q}$
 &$d=0,q=3,A_q=0$ only \\
Triangle-free $F$, non-discrete $d\ge1$, path scale
 &$N^{1-d}\varepsilon^{q-3+A_q}$
 &$d=1,q=3,A_q=0$ only; all others $O(N^{-1}+\varepsilon)$ \\
Girth-$g$ $F$, discrete quotient, shortest-cycle scale
 &$\varepsilon^{q-g+A_q}$
 &$q=g,A_q=0$ only; all others $O(\varepsilon)$ \\
Triangle-free $F$, non-discrete $d\ge1$, four-cycle scale
 &$N^{1-d}\varepsilon^{q-3+A_q}/u$
 &$d=1,q=3,A_q=0$ can only be rooted $K_3$, which is constant when $\tau_R$ is constant \\
\bottomrule
\end{tabular}
\end{center}
The right-hand sides are obtained by division of
\eqref{eq:quotient} by the respective scales; use $u/N\le\varepsilon$
and $u\ge1$. The discrete support is excluded from the second row because at
larger girth it may contribute at a different principal scale. If $A_q=0$, then $m=q+c-1$; every additional nontrivial
component has at least two vertices, forcing $c=1$, $m=q$ and a
chordless rooted cycle. For $d=1$ and $q=3$, there are exactly four
vertices in the support and the only triangle-producing independent
identification in a triangle-free support is the endpoint pair of
$P_4$. Hence the table exhausts the zero-loss possibilities without
assuming a bound on the order of $F$.

The coefficients require counting the \emph{original} source roots,
not only the quotient vertices. The following ledger also records the
sign from the chosen support and its partition:
\begin{center}
\small
\begin{tabular}{@{}lcccc@{}}
\toprule
Type & Source-root weight & $(-1)^{|E'|}$ & $\mu(\pi)$ & Varying rooted term \\
\midrule
$K_3$ & $3$ & $-1$ & $1$ & $2\Delta\tau_R$ \\
Endpoint-identified $P_4$ & $2+1+1=4$ & $-1$ & $-1$ & $2\Delta\tau_R$ \\
$C_4$ & $4$ & $+1$ & $1$ & $2\Delta c_4^R$ \\
$C_g$ & $g$ & $(-1)^g$ & $1$ & $\Delta z_g^R$ \\
\bottomrule
\end{tabular}
\end{center}
Multiplying each row gives respectively $-6c_3(F)$,
$8p_4(F)$, $8c_4(F)$ and $(-1)^ggc_g(F)$ after summing copies
of the corresponding support in $F$. For the $C_4$ row,
$z_4=2c_4+2\omega^2-\omega$; the regular-degree terms cancel in
root differences. The $P_4$ contraction uses the four \emph{source}
root positions, although its quotient triangle has only three
vertices. This bookkeeping is independent of the finite seed certificates.

\Needspace{18\baselineskip}
\section{Finite seeds from a block substitution}\label{app:embedded}
The triangle alphabet is obtained from a threshold graph by deleting one vertex and repairing its neighborhood with a matching. A block substitution then produces the four-cycle alphabet from the same complementary pair. The substitution is specified by the five-by-five matrix below. All indices in this appendix start at zero.

\subsection*{The triangle alphabet}
Put $m=63$ and define two $m\times m$ matrices by
\[
 C_{ij}=\mathbf1_{\{i\ne j,\ i+j<62\}},\qquad
 X_{ij}=\mathbf1_{\{i+j>62\ \text{or}\ i=j<31\}}.
\]
In $X$, the notation $(a,b;c,d)$ denotes a switch from ones at $(a,c),(b,d)$ to ones at $(a,d),(b,c)$, with the other two entries initially zero. Perform $(a,a+2;a,62-a)$ for $a=0,\ldots,28$, in this order, followed by $(2,30;0,30)$ and $(16,48;14,31)$. These are the threshold-block switches used in~\cite{ZhangK3}; they are specified here in full. Each preserves the row and column sums of $X$. The graph with adjacency matrix
\[
 G_* = \begin{pmatrix} C&X\\X^T&C\end{pmatrix}
\]
is therefore $62$-regular on the labels $0,\ldots,125$.

Delete vertex $42$, whose neighborhood is $\{0,\ldots,19\}\cup\{84,\ldots,125\}$. To specify a matching on this set, put
\[
 L_i=\begin{cases}i,&i<20,\\84+i,&20\le i<31,\end{cases}
 \qquad R_j=\begin{cases}84+j,&j<20,\\95+j,&20\le j<31,\end{cases}
\]
and define the permutations
\[
 f(i)=\begin{cases}i,&i<20,\\20+\bigl((9(i-20)+7)\bmod11\bigr),&20\le i<31,\end{cases}
\]
\[
 \sigma=(0\ 26\ 6)(2\ 11)(3\ 5)(8\ 13)(9\ 19)(24\ 29).
\]
Each cycle is read from left to right. Add the edges
\begin{equation}\label{eq:compactmatching}
 \{L_i,R_{f(\sigma(i))}\},\qquad 0\le i<31.
\end{equation}
The lists $L$ and $R$ partition the deleted neighborhood, and $f\circ\sigma$ is a permutation, so these edges form a perfect matching. The threshold rule and the displayed switches show that each pair was a nonedge. Every former neighbor loses one edge and gains one edge. Renumbering the remaining vertices in increasing order yields a $125$-vertex $62$-regular graph $A$; put $B=\overline A$.

For any $d$-regular graph $H$ on $2d+1$ vertices,
\[
 \tau_H(v)+\tau_{\overline H}(v)=\frac{d(d-2)}2.
\]
Indeed, writing $S=N_H(v)$ and $T=V(H)\setminus(S\cup\{v\})$, the degree sums give $e_H(T)-e_H(S)=d/2$, and $\tau_{\overline H}(v)=\binom d2-e_H(T)$. Thus $\tau_A(v)+\tau_B(v)=1860$. The exact finite check supplied with the source establishes the joint injectivity and range of these two triangle alphabets.

\subsection*{The block substitution}
Let
\begin{equation}\label{eq:controlfive}
 T=\begin{pmatrix}
 1&0&0&0&0\\
 1&1&0&0&0\\
 1&0&1&0&0\\
 0&1&1&1&0\\
 0&1&0&1&1
 \end{pmatrix},\qquad
 M_U=T\otimes A+(J_5-T)\otimes B,\quad
 M_V=T\otimes B+(J_5-T)\otimes A.
\end{equation}
Here $A,B$ also denote adjacency matrices. The two bipartite graphs with incidence matrices $M_U,M_V$ are denoted by $U,V$. Each row and column contains five blocks of degree $62$, so both graphs are $310$-regular on $1250$ vertices. They have disjoint edge sets because $A$ and $B$ do. Their common bipartition has $625$ vertices in each class. No regularity condition is needed on $T$ itself.

The asymmetry of $T$ is essential to this method of separating both bipartition classes: if $T$ were symmetric, then $M_U$ would be symmetric, giving equal rooted four-cycle counts at corresponding left and right vertices. The precise injectivity of the displayed choice is verified below.

We first establish the block fourth-moment identity for an arbitrary control size $q$. For the present $q=5$ substitution, its seed-dependent terms reduce to the rooted triangle count and one fourth diagonal moment of $A$.
\begin{lemma}[Block fourth-moment identity]\label{lem:blockfourth}
Let $A,B$ be adjacency matrices of complementary $d$-regular graphs on $n=2d+1$ vertices, and let $Z$ be any $q\times q$ zero--one matrix. Put $E=A-B$ and
\[
 \mathsf M=J_q\otimes B+Z\otimes E,\qquad
 P_0=qJ_q,\quad P_1=ZJ_q+J_qZ^T,\quad P_2=ZZ^T.
\]
For $j=0,\ldots,4$, let $a_j(i)$ be the $i$th diagonal entry of the corresponding matrix in
\[
\begin{gathered}
 P_0^2,\quad P_0P_1+P_1P_0,\quad P_1^2+P_0P_2+P_2P_0,\\
 P_1P_2+P_2P_1,\quad P_2^2.
\end{gathered}
\]
In the bipartite graph with incidence matrix $\mathsf M$, the rooted four-cycle count at the left vertex $(i,x)$ is
\begin{equation}\label{eq:blockcount}
 c_4(i,x)=\frac12\left(\sum_{j=0}^4 a_j(i)(B^{4-j}E^j)_{xx}
                       -2(qd)^2+qd\right).
\end{equation}
Right-vertex counts are obtained by replacing $Z$ with $Z^T$. For the complementary block substitution, replace $Z$ with $J_q-Z$.
\end{lemma}
\begin{proof}
Regularity gives $AJ_n=J_nA=dJ_n$, so $A$ and $B=J_n-I_n-A$ commute. Hence $B$ and $E$ commute, and
\[
 \mathsf M\mathsf M^T=P_0\otimes B^2+P_1\otimes BE+P_2\otimes E^2.
\]
Squaring this identity yields the five coefficient matrices above. Put $H=\mathsf M\mathsf M^T$ and $\delta=qd$; for each left vertex $v$ we have $H_{vv}=\delta$ and $\sum_w H_{vw}=\delta^2$. Counting a four-cycle by the vertex opposite $v$ gives
\[
 c_4(v)=\sum_{w\ne v}\binom{H_{vw}}2
       =\frac{(H^2)_{vv}-2\delta^2+\delta}{2},
\]
which proves~\eqref{eq:blockcount}. Since $A$ and $B$ are symmetric, $\mathsf M^T$ is the same substitution with $Z^T$. Interchanging $A$ and $B$ replaces $Z$ by $J_q-Z$.
\end{proof}

For $n=125$ and $d=62$, write $t_x=\tau_A(x)$ and $w_x=(A^4)_{xx}$. Direct expansion of $B=J-I-A$ and $E=2A+I-J$ gives
\[
\begin{array}{c|l}
 j&(B^{4-j}E^j)_{xx}\\ \hline
 0&w_x+8t_x-7440\\
 1&-2w_x-14t_x+249488\\
 2&4w_x+24t_x-495318\\
 3&-8w_x-40t_x+983258\\
 4&16w_x+64t_x-1951636
\end{array}
\]
Thus the joint four-cycle counts are twenty explicit affine functions of $(t_x,w_x)$, evaluated at the $125$ seed vertices. The required injectivity concerns this specific pair $A,T$; the supplementary verifier certifies it using both Eq.~\eqref{eq:blockcount} and direct common-neighbor counts in $U,V$.

The resulting exact rooted-count ranges are recorded below.
\begin{center}
\small
\begin{tabular}{@{}lrrr@{}}
\toprule
Rooted statistic & Distinct values & Minimum & Maximum\\
\midrule
$\tau_A\cup\tau_B$ &250&464&1396\\
$\tau_{H_0}$ &31250&7971214&22534646\\
$c_4^U\cup c_4^V$ &2500&7783825&8618509\\
\bottomrule
\end{tabular}
\end{center}
The $H_0$ values follow from Eq.~\eqref{eq:trianglelex}. The remaining injectivity assertions are certified by the exact-integer verifier \texttt{verify\_general\_seeds.py} in the supplementary archive. It reconstructs the seeds from the displayed rules and compares Eq.~\eqref{eq:blockcount} with direct common-neighbor enumeration; all 2500 rooted four-cycle counts must agree entrywise and be distinct. Run with Python~3.10 or later, without optimization flags so assertions remain active. The program is reproduced in the PDF only when the optional finite-code switch is enabled.

\Needspace{20\baselineskip}
\section{Self-contained 29-row finite certificate}\label{app:qt-certificate}
The following 29 hexadecimal integers are ordered by vertex labels $0,1,\ldots,28$. An integer's bit $j$ is one precisely when its row-vertex is adjacent to vertex $j$. All the numbers fit within 29 bits.
\begin{center}
\begin{tabular}{r@{\quad}l@{\qquad}r@{\quad}l@{\qquad}r@{\quad}l}
\toprule
$i$& mask&$i$&mask&$i$&mask\\\midrule
0&\texttt{6c622f6}&10&\texttt{1d1f4b40}&20&\texttt{118fde08}\\
1&\texttt{6c023fd}&11&\texttt{1d1f5460}&21&\texttt{1b0ff120}\\
2&\texttt{2c033fb}&12&\texttt{1df9090c}&22&\texttt{282b3d7}\\
3&\texttt{49453f6}&13&\texttt{12e103b7}&23&\texttt{25031ff}\\
4&\texttt{2c423ef}&14&\texttt{1d3d8d08}&24&\texttt{1c3bdc40}\\
5&\texttt{ea0a89f}&15&\texttt{d7f42a0}&25&\texttt{e223f7}\\
6&\texttt{7c00e9f}&16&\texttt{193efc00}&26&\texttt{908dd6b}\\
7&\texttt{2c0a37f}&17&\texttt{1b7d8c01}&27&\texttt{152fdc20}\\
8&\texttt{6e0749e}&18&\texttt{183bcc19}&28&\texttt{93f7e00}\\
9&\texttt{1250a4df}&19&\texttt{1d37dc00}&&\\\bottomrule
\end{tabular}
\end{center}
The supplementary script \texttt{verify\_triangle29.py} enumerates all unordered vertex triples from these masks and checks symmetry, zero diagonal, 14-regularity and the exact count sets in Lemma~\ref{lem:qt-seed}. It requires no external data.

\Needspace{8\baselineskip}
\section*{Statements and Declarations}
\noindent\textbf{Funding.} The author received no external funding for this research.

\noindent\textbf{Competing interests.} The author declares no financial or non-financial interests directly related to this work.

\noindent\textbf{Data and code availability.} The finite seed specifications appear in Appendices~\ref{app:embedded} and~\ref{app:qt-certificate}. Deterministic Python verification scripts are supplied as supplementary source files and require no external datasets.

\noindent\textbf{Use of generative AI tools.} ChatGPT was used as a supporting tool for mathematical exploration, manuscript preparation, and verification-code assistance. The author reviewed the mathematical arguments and takes full responsibility for the results and final manuscript.

\Needspace{9\baselineskip}
\medskip
\begingroup\fontsize{9.5}{10.5}\selectfont\raggedright
\setlength{\bibsep}{0pt}

\endgroup
\end{document}